\documentclass[12pt]{amsart} 

\makeatletter
\newcommand\mathcircled[1]{%
	\mathpalette\@mathcircled{#1}%
}
\newcommand\@mathcircled[2]{%
	\tikz[baseline=(math.base)] \node[draw,ellipse,inner sep=1pt] (math) {$\m@th#1#2$};%
}

\usepackage{multicol}
\usepackage{latexsym}
\usepackage{psfrag}
\usepackage{amsmath}
\usepackage{amssymb}
\usepackage{epsfig}
\usepackage{amsfonts}
\usepackage{amscd}
\usepackage{mathrsfs}
\usepackage{graphicx}
\usepackage{enumerate}
\usepackage[autostyle=false, style=english]{csquotes}
\MakeOuterQuote{"}
\usepackage{ragged2e}
\usepackage[all]{xy}
\usepackage{mathtools}

\newlength\ubwidth

\usepackage[dvipsnames]{xcolor}

\usepackage{placeins}
\usepackage{tikz}
\usetikzlibrary{shapes, positioning}
\usetikzlibrary{matrix,shapes.geometric,calc,backgrounds}
\usetikzlibrary{arrows.meta} 
\usetikzlibrary{shapes.geometric, arrows, positioning, decorations.pathreplacing} 

\newtheorem{theorem}{Theorem}

\usepackage{amsmath,amssymb,amsthm}
\usepackage[dvipsnames]{xcolor}
\usepackage[colorlinks=true,citecolor=RoyalBlue,linkcolor=red,breaklinks=true]{hyperref}

\usepackage{graphicx}
\usepackage{mathdots} 
\usepackage[margin=2.5cm]{geometry}
\usepackage{ytableau}

\usepackage{cancel}

\usepackage{blkarray}
\usepackage{booktabs}

\dedicatory{Dedicated to George E. Andrews and Bruce C. Berndt 
for their 85th birthday.}

\begin{document}

\title[Finding systems of equations]
{Finding systems of functional equations for Andrews-Gordon type series}

\author[K{\i}l{\i}\c{c}]{Yal\c{c}{\i}n Can K{\i}l{\i}\c{c}}
\address{Yal\c{c}{\i}n Can K{\i}l{\i}\c{c}, Faculty of Engineering and Natural Sciences, 
    Sabanc{\i} University, Tuzla, Istanbul 34956, Turkey}
\email{yalcinkilic@sabanciuniv.edu}

\author[Kur\c{s}ung\"{o}z]{Ka\u{g}an Kur\c{s}ung\"{o}z}
\address{Ka\u{g}an Kur\c{s}ung\"{o}z, Faculty of Engineering and Natural Sciences, 
    Sabanc{\i} University, Tuzla, Istanbul 34956, Turkey}
\email{kursungoz@sabanciuniv.edu}


      
\thanks{This research was partially supported by T\"{U}B\.{I}TAK grant no. 122F136. }

\date{2025}

\begin{abstract}
 We develop a search algorithm for systems of $q$-difference equations 
 satisfied by Andrews-Gordon type double series.  
 We then couple the search algorithm with Euler's algorithm 
 for finding infinite products to narrow the search space.  
 We exemplify some findings of the algorithm, 
 along with their proofs.  
 We also explain some of the double series in a base partition and moves framework.  
\end{abstract}

\maketitle

\section{Introduction and statement of results}
\label{secIntro}

A.O.L. Atkin is quoted as saying 
``\ldots it is often more difficult to discover results 
in this subject than to prove them \ldots''~\cite{G_Atkin_quote}.  
After decades of supporting evidence, 
largely thanks to George Andrews (e.g.~\cite{A_Capp}, to cite but one significant example), 
the situation began to reverse.  
Inspired by Capparelli's work~\cite{Capp_th}, Kanade and Russell 
came up with a systematic way to discover partition identities in 2015~\cite{KR_15}.  
Their mod 12 conjectures have been proven~\cite{BJM, Ros}, 
but the mod 9 conjectures are still open.  

The main aim of this paper is to treat the series 
\begin{align}
\label{srMainDouble}
 S_{C_1, C_2}(x) = S_{C_1, C_2}(x; q) = \sum_{m, n \geq 0} 
  \frac{ q^{ B_{11} \binom{m+1}{2} + B_{22} \binom{n+1}{2} + B_{12}mn + C_1 m + C_2 n } 
    x^{ D_1 m + D_2 n } }
   { (q^{K_1}; q^{K_1})_m (q^{K_2}; q^{K_2})_n }
\end{align}
with computer algebra, 
and discover systems of functional equations 
that uniquely determine $S_{C_1, C_2}(x)$ 
for arbitrary but fixed parameters 
$B_{11}$, $B_{22}$, $B_{12}$, $C_1$, $C_2$, $D_1$, $D_2$, $K_1$ and $K_2$, 
along with auxiliary series of the same type.  
Here, 
\begin{align}
\nonumber 
  (a; q)_n = \prod_{j = 1}^n ( 1 - a q^{j-1})
\end{align}
is the $q$-Pochhammer symbol.  
Later in the paper, we will use 
\begin{align}
\nonumber 
  (a; q)_\infty = \lim_{n \to \infty} (a; q)_n.  
\end{align}
It suffices to take $\vert q \vert < 1$ so that 
all single or multiple series and the infinite products converge absolutely~\cite{GR}.  
The convergence condition is mentioned for the sake of completeness, 
as the results in this paper do not require convergence per se.  
They can be proven formally~\cite{Wilf_GFunc}.  
The parameters $D_1$, $D_2$, $K_1$ and $K_2$ are positive integers. 
The parameters $B_{11}$, $B_{22}$, $B_{12}$, $C_1$, and $C_2$ are chosen such that 
$S_{C_1, C_2}(1)$ is a power series in $q$ such that $S_{C_1, C_2}(0) = 1$.  

The secondary aim of the paper is to give examples 
of the findings of the search algorithm along with their proofs.  
The applications are interpretations of instances of \eqref{srMainDouble} 
as partition generating functions, 
as well as sum-product identities that can be interpreted as partition identities.  

An integer partition of a non-negative integer $n$ 
is an unordered sum of positive integers that add up to $n$.  
For instance, 4 has the five partitions 
\begin{align}
\nonumber 
  1+1+1+1, \quad 1+1+2, \quad 2+2, \quad 1+3, \quad 4.  
\end{align}
We agree that zero has one partition, namely empty partition~\cite{A_thebluebook}.  
Depending on the context, which will be made clear, 
we will allow a limited number of zeros to appear as summands, such as $0+0+1+1+2$.  
This is to fix the number of parts (also called the \emph{length}).  
$1+1+2$ has length three, where $0+0+1+1+2$ has length five.  
It may also be convenient to 
write partitions as a list, 
or drop the plus signs in between and 
represent the repeated parts stacked 
on top of each other, such as $(1,1,2)$ or 
$\begin{array}{cc} 1 & \\ 1 & 2  \end{array}$
instead of $1+1+2$.  

The result samples are as follows.  

\begin{theorem}
\label{thmFirstApplication}
 Let $f(m,n)$ be partitions of $n$ into $m$ parts 
 which do not repeat more than twice, 
 such that part sizes are at least two apart, 
 and repeating part sizes are at least three apart.  
 Then, 
 \begin{align}
 \label{eqFirstApplication}
  \sum_{m, n \geq 0} f(m, n) x^m q^n 
  = \sum_{m, n} \frac{ q^{ 6\binom{m+1}{2} + 2\binom{n+1}{2} + 2mn - 4m - n } x^{2m+n} }
    { ( q^2; q^2)_m ( q; q)_n }
 \end{align}
\end{theorem}

\begin{theorem}
\label{thmYalcinMain}
 For any non-negative integer $n$, the following sets have the same cardinalities.  
 \begin{enumerate}[1.]
  \item The set of ordinary partitions of $n$ in which each part appears at most three times.  
  \item The set of bicolored partitions of $n$ such that 
    \begin{enumerate}[(a)]
      \item Neither blue nor red parts repeat.   
      \item When we go through parts from the smallest to the largest, 
      for each red $\textcolor{red}{b}$, there must be a blue $\textcolor{blue}{b}$ 
      or a blue $\textcolor{blue}{(b+1)}$.  
      Moreover, these blue parts should be different for each red $\textcolor{red}{b}$.  
    \end{enumerate}
  \item The set of bicolored partitions of $n$ such that 
    \begin{enumerate}[(a)]
      \item Both blue and red parts appear at most once.  
      \item Red $\textcolor{red}{1}$ does not appear as a part.  
      \item For any blue $\textcolor{blue}{b}$, we cannot have a red $\textcolor{red}{(b+1)}$ 
        or a red $\textcolor{red}{(b+2)}$.  
      \item In any sublist of successive parts ending with a blue $\textcolor{blue}{b}$, 
        we cannot have differences $[1,2^{\ast}]$ with successive red parts.  
        Thus, we cannot have a red $\textcolor{red}{(b-1)} $ and a blue $ \textcolor{blue}{b}$,
        or a red $\textcolor{red}{(b-3)}$ a red $\textcolor{red}{(b-2)}$ and a blue $\textcolor{blue}{b}$,
        or a red $\textcolor{red}{(b-5)}$ a red $\textcolor{red}{(b-4)}$ a red
        $\textcolor{red}{(b-2)}$ and a blue $\textcolor{blue}{b}$ etc.
    \end{enumerate}
 \end{enumerate}
\end{theorem}
Here, the notation $2^{\ast}$ means any number of appearances of $2$, including zero.  

{\bf Example:}
Let us give some examples and non-examples to clarify the conditions in Theorem \ref{thmYalcinMain}:

\begin{enumerate}
\item $\textcolor{red}{1} + \textcolor{red}{2} + \textcolor{blue}{2} $ does not satisfy the second condition of the second partition class, since for the red $\textcolor{red}{1}$ and for the red $\textcolor{red}{2}$ we have the blue $\textcolor{blue}{2}$, however there should be a different blue part for both of them.
\item $\textcolor{red}{1} + \textcolor{red}{2} + \textcolor{blue}{2} + \textcolor{blue}{3}$ satisfies the conditions of the second partition class.
\item $\textcolor{red}{3} + \textcolor{blue}{5}$ satisfies the conditions of the third partition class.
\item $\textcolor{red}{2} + \textcolor{red}{3} + \textcolor{blue}{5}$ does not satisfy the conditions of the third partition class, since it violates the last condition. 
\end{enumerate}

\begin{theorem}
\label{thmAltRRG3} 
 {\allowdisplaybreaks
 \begin{align}
 \nonumber
  \sum_{m, n \geq 0} \frac{ (-1)^n q^{ 2 \binom{m+1}{2} + 2 \binom{n+1}{2} + 2mn - m - n } }
    { (q; q)_m (q^2; q^2)_n }
  & = \frac{ 1 }{ ( q^2, q^4, q^{10}, q^{12}; q^{14} )_\infty } \\ 
 \nonumber
  \sum_{m, n \geq 0} \frac{ (-1)^n q^{ 2 \binom{m+1}{2} + 2 \binom{n+1}{2} + 2mn + n } }
    { (q; q)_m (q^2; q^2)_n }
  & = \frac{ 1 }{ ( q^2, q^6, q^8, q^{12}; q^{14} )_\infty } \\ 
 \nonumber
  \sum_{m, n \geq 0} \frac{ (-1)^n q^{ 2 \binom{m+1}{2} + 2 \binom{n+1}{2} + 2mn + m + n } }
    { (q; q)_m (q^2; q^2)_n }
  & = \frac{ 1 }{ ( q^4, q^6, q^8, q^{10}; q^{14} )_\infty } 
 \end{align}} 
\end{theorem}

The right-hand sides are the $q \leftarrow q^2$ dilation of the Rogers-Ramanujan-Gordon 
identities for $k = 3$~\cite{G}.  
We must impress that 
it is not obvious that the left-hand sides are even functions of $q$.  

\begin{theorem}
\label{thmAltAlladiGordonDilation}
 \begin{align}
 \nonumber 
  \sum_{m, n \geq 0} \frac{ (-1)^n q^{ 9 \binom{m+1}{2} + 6 \binom{n+1}{2} + 6mn - n } }
    { (q^3; q^3)_m (q^2; q^2)_n }
  = \left( q, q^5; q^6 \right)_\infty
 \end{align}
\end{theorem}

The rest of the paper is organized as follows.  
In Section \ref{secSrchAlgo}, the search algorithm is explained with a detailed example.  
In Section \ref{secFirstApplication}, 
Theorem \ref{thmFirstApplication} is proven, 
and it is also explained using a base partition and moves framework.  
The search space is narrowed thanks to Euler's algorithm~\cite{A_anotherbook, KR_15}
for finding infinite products, hence partition identities, 
in Section \ref{secNarrowingSrchSpc}.  
Two outputs of the enhanced search algorithm are given, 
and one series is explained using a base partition and moves framework.  
Then, Theorem \ref{thmYalcinMain} is proven.  
A low hanging fruit is collected in Section \ref{secAltSeries}, 
and series that are not evidently positive are considered.  
Combinatorial proofs 
of Theorems \ref{thmAltRRG3} and \ref{thmAltAlladiGordonDilation} are given.  
We conclude with possible directions of future research in Section \ref{secFutureWork}.  

The idea elaborated on in this paper is not exactly new.  
It was used in different forms in~\cite{Chern} and in~\cite{KR_20}, as notable examples.  
We will mention the connections and differences in Section \ref{secSrchAlgo}.

\section{The search algorithm}
\label{secSrchAlgo}

The double series \eqref{srMainDouble} satisfies the following 
\emph{primary} $q$-contiguous equations 
for an arbitrary but fixed positive integer $\gamma$. 
\begin{align}
\nonumber 
 S_{C_1, C_2}(x) - S_{C_1 + K_1, C_2}(x) 
 & = x^{D_1} q^{ B_{11} + C_1 } 
  S_{C_1 + B_{11} - \gamma D_1, C_2 + B_{12} - \gamma D_2 }(x q^{\gamma}) \\
\nonumber 
 S_{C_1, C_2}(x) - S_{C_1 , C_2 + K_2}(x) 
 & = x^{D_2} q^{ B_{22} + C_2 } 
  S_{C_1 + B_{12} - \gamma D_1, C_2 + B_{22} - \gamma D_2 }(x q^{\gamma}) \\
\label{fneqPrimary}
  S_{C_1, C_2}(x) & = S_{C_1 - \gamma D_1 , C_2 - \gamma D_1}(x q^{\gamma}) 
\end{align}
We justify the first one. 
{\allowdisplaybreaks
\begin{align}
\nonumber 
 & S_{C_1, C_2}(x) - S_{C_1 + K_1, C_2}(x) = \sum_{m, n \geq 0} 
  \frac{ q^{ B_{11} \binom{m+1}{2} + B_{22} \binom{n+1}{2} + B_{12}mn + C_1 m + C_2 n } 
    x^{ D_1 m + D_2 n } \left( 1 - q^{K_1 m} \right) }
   { (q^{K_1}; q^{K_1})_m (q^{K_2}; q^{K_2})_n } \\
\nonumber 
 & = \sum_{m \geq 1, n \geq 0} 
  \frac{ q^{ B_{11} \binom{m+1}{2} + B_{22} \binom{n+1}{2} + B_{12}mn + C_1 m + C_2 n } 
    x^{ D_1 m + D_2 n } }
   { (q^{K_1}; q^{K_1})_{m - 1} (q^{K_2}; q^{K_2})_n } \\
\nonumber 
 & = \sum_{m, n \geq 0} 
  \frac{ q^{ B_{11} \binom{m+2}{2} + B_{22} \binom{n+1}{2} + B_{12}(m+1)n + C_1 (m+1) + C_2 n } 
    x^{ D_1 ( m + 1) + D_2 n } }
   { (q^{K_1}; q^{K_1})_{m} (q^{K_2}; q^{K_2})_n }  \\
\nonumber 
 & = \sum_{m, n \geq 0} 
  \frac{ q^{ B_{11} \binom{m+1}{2} + B_{11}(m+1) + B_{11} + B_{22} \binom{n+1}{2} + B_{12}mn 
    + B_{12}n + C_1 m + C_1 + C_2 n } 
    x^{ D_1 m + D_1 + D_2 n } }
   { (q^{K_1}; q^{K_1})_{m} (q^{K_2}; q^{K_2})_n } \\
\nonumber 
 & = x^{D_1} q^{ C_1 + B_{11} } \sum_{m, n \geq 0} 
  \frac{ q^{ B_{11} \binom{m+1}{2} + B_{22} \binom{n+1}{2} + B_{12}mn 
    + (C_1 + B_{11} )m + (C_2 + B_{12}) n } 
    x^{ D_1 m + D_2 n } }
   { (q^{K_1}; q^{K_1})_{m} (q^{K_2}; q^{K_2})_n } \\
\nonumber 
 & = x^{D_1} q^{ C_1 + B_{11} } \sum_{m, n \geq 0} 
  \frac{ q^{ B_{11} \binom{m+1}{2} + B_{22} \binom{n+1}{2} + B_{12}mn 
    + (C_1 + B_{11} -\gamma D_1 )m + (C_2 + B_{12} - \gamma D_2) n } 
    (x q^{ \gamma } )^{ D_1 m + D_2 n } }
   { (q^{K_1}; q^{K_1})_{m} (q^{K_2}; q^{K_2})_n } \\
\nonumber 
 & = x^{D_1} q^{ B_{11} + C_1 } 
  S_{C_1 + B_{11} - \gamma D_1, C_2 + B_{12} - \gamma D_2 }(x q^{\gamma})
\end{align} }
The second one follows by symmetry, 
and the third one is but the final twist in the proof of the first one.  
The computations above make clear that the parameters 
$B_{11}$, $B_{12}$, $B_{22}$, $D_1$, $D_2$, $K_1$ and $K_2$ do not change. 
This is why only $C_1$ and $C_2$ are highlighted in \eqref{srMainDouble}. 

The first indices appearing in the functional equations \eqref{fneqPrimary} are 
\begin{align}
\label{listIndices1}
  C_1, \quad C_1 + K_1, \quad 
  C_1 + B_{11} - \gamma D_1, \quad 
  C_1 + B_{12} - \gamma D_1, \quad 
  C_1 - \gamma D_1, 
\end{align}
and the second indices are 
\begin{align}
\label{listIndices2}
  C_2, \quad C_2 + K_2, \quad 
  C_2 + B_{12} - \gamma D_2, \quad 
  C_2 + B_{22} - \gamma D_2, \quad 
  C_2 - \gamma D_2.  
\end{align}
We choose and fix intervals for these indices, 
say, $[m_1, M_1]$ for the indices in \eqref{listIndices1}, 
and $[m_2, M_2]$ for the indices in \eqref{listIndices2}.  
Then, we list all instances of the functional equations in \eqref{fneqPrimary} 
all of the first (respectively, the second) indices of which 
fall in $[m_1, M_1]$ (respectively, in $[m_2, M_2]$).  
It is best to work on a running example.  
It is a heavily guided and well studied one 
for demonstration purposes~\cite{A_PNAS, A_thebluebook}.  

We know that the three series in Andrews-Gordon identities 
for $k = 3$~\cite{A_PNAS} are 
\begin{align}
\nonumber 
  \sum_{m, n \geq 0} \frac{ q^{ (m+n)^2 + m^2 } x^{2m + n} }{ (q; q)_m (q; q)_n }
  = \sum_{m, n \geq 0} \frac{ q^{ 4 \binom{m+1}{2} + 2 \binom{n+1}{2} + 2 mn - 2m -n} x^{2m + n} }
    { (q; q)_m (q; q)_n }
  = S_{-2, -1}(x), 
\end{align}
\begin{align}
\nonumber 
  \sum_{m, n \geq 0} \frac{ q^{ (m+n)^2 + m^2 + m } x^{2m + n} }{ (q; q)_m (q; q)_n }
  = \sum_{m, n \geq 0} \frac{ q^{ 4 \binom{m+1}{2} + 2 \binom{n+1}{2} + 2 mn - m -n} x^{2m + n} }
    { (q; q)_m (q; q)_n }
  = S_{-1, -1}(x), 
\end{align}
\begin{align}
\nonumber 
  \sum_{m, n \geq 0} \frac{ q^{ (m+n)^2 + m^2 + 2m + n} x^{2m + n} }{ (q; q)_m (q; q)_n }
  = \sum_{m, n \geq 0} \frac{ q^{ 4 \binom{m+1}{2} + 2 \binom{n+1}{2} + 2 mn } x^{2m + n} }
    { (q; q)_m (q; q)_n }
  = S_{0, 0}(x).  
\end{align}
The common parameters are $(B_{11}, B_{22}, B_{12}, D_1, D_2, K_1, K_2, \gamma)$ 
$=(4,2,2,2,1,1,1,1)$. 
Andrews~\cite{A_PNAS, A_thebluebook} showed that
\begin{align}
\nonumber
  S_{-2, -1}(x) - S_{-1, -1}(x) & = x^2 q^2 S_{0, 0}(xq), \\
\nonumber
  S_{-1, -1}(x) - S_{0, 0}(x) & = x q S_{-1, -1}(xq), \\
\label{fneq_A_PNAS}
  S_{0, 0}(x) & = S_{-2, -1}(xq).
\end{align}
The first and the last of these equations fit \eqref{fneqPrimary},
but the second one does not.
Andrews uses an inductive argument to demonstrate the second equation in \eqref{fneq_A_PNAS}.
The algorithm we outline below, with all due respect, bypasses ingenuity.
It tries and finds a way to express the second equation in \eqref{fneq_A_PNAS}
as a linear combination of instances of \eqref{fneqPrimary}.

We choose the interval $[-2, 1]$ for the first index in \eqref{srMainDouble}, 
and $[-1, 1]$ for the second.  
All instances of \eqref{fneqPrimary} whose indices 
fall in the respective intervals are as follows.  
{\allowdisplaybreaks
\begin{align}
\nonumber 
 & S_{ -2, -1}(x) - S_{ -1, -1}(x) - x^{ 2 } q^{ 2 } S_{ 0, 0}(xq) = 0 \\
\nonumber 
 & S_{ -2, 0}(x) - S_{ -1, 0}(x) - x^{ 2 } q^{ 2 } S_{ 0, 1}(xq) = 0 \\
\nonumber 
 & S_{ -1, -1}(x) - S_{ 0, -1}(x) - x^{ 2 } q^{ 3 } S_{ 1, 0}(xq) = 0 \\
\nonumber 
 & S_{ -1, 0}(x) - S_{ 0, 0}(x) - x^{ 2 } q^{ 3 } S_{ 1, 1}(xq) = 0 \\[2mm] 
\nonumber 
 & S_{ -2, -1}(x) - S_{ -2, 0}(x) - x q S_{ -2, 0}(xq) = 0 \\
\nonumber 
 & S_{ -2, 0}(x) - S_{ -2, 1}(x) - x q^{ 2 } S_{ -2, 1}(xq) = 0 \\
\nonumber 
 & S_{ -1, -1}(x) - S_{ -1, 0}(x) - x q S_{ -1, 0}(xq) = 0 \\
\nonumber 
 & S_{ -1, 0}(x) - S_{ -1, 1}(x) - x q^{ 2 } S_{ -1, 1}(xq) = 0 \\
\nonumber 
 & S_{ 0, -1}(x) - S_{ 0, 0}(x) - x q S_{ 0, 0}(xq) = 0 \\
\nonumber 
 & S_{ 0, 0}(x) - S_{ 0, 1}(x) - x q^{ 2 } S_{ 0, 1}(xq) = 0 \\
\nonumber 
 & S_{ 1, -1}(x) - S_{ 1, 0}(x) - x q S_{ 1, 0}(xq) = 0 \\
\nonumber 
 & S_{ 1, 0}(x) - S_{ 1, 1}(x) - x q^{ 2 } S_{ 1, 1}(xq) = 0 \\[2mm] 
\nonumber 
 & S_{ 0, 0}(x) - S_{ -2, -1}(xq) = 0 \\
\nonumber 
 & S_{ 0, 1}(x) - S_{ -2, 0}(xq) = 0 \\
\nonumber 
 & S_{ 1, 0}(x) - S_{ -1, -1}(xq) = 0 \\
\label{fneqRunningExample} 
 & S_{ 1, 1}(x) - S_{ -1, 0}(xq) = 0 
\end{align}}
16 functional equations are listed for a total of 24 series, namely 
\begin{align}
\nonumber 
  S_{ -2, -1}(x), \; 
  S_{ -1, -1}(x), \; 
  S_{ -2, 0}(x), \; 
  S_{ 0, -1}(x), \; 
  S_{ -1, 0}(x), \; 
  S_{ -2, 1}(x), \; 
  S_{ 1, -1}(x), \; 
  S_{ 0, 0}(x),
\end{align}
\begin{align}
\nonumber 
  S_{ -1, 1}(x), \; 
  S_{ 1, 0}(x), \; 
  S_{ 0, 1}(x), \; 
  S_{ 1, 1}(x), 
\end{align}
and the same series with $x$ replaced by $xq$.  
We form linear combinations of the 16 equations in \eqref{fneqRunningExample} 
by multiplying the equations by $t_1$, $t_2$, \ldots, $t_{16}$ 
in their respective order, 
add all side by side, 
and collect coefficients of the like series, 
to get the following.  
{\allowdisplaybreaks
\begin{align}
\nonumber
  & \Big( t_1 + t_5 \Big) S_{ -2, -1}(x)
  + \Big( -t_1 + t_3 + t_7 \Big) S_{ -1, -1}(x)
  + \Big( t_2 - t_5 + t_6 \Big) S_{ -2, 0}(x) \\
\nonumber
  & + \Big( -t_3 + t_9 \Big) S_{ 0, -1}(x)
  + \Big( - t_2 + t_4 - t_7 + t_8 \Big) S_{ -1, 0}(x)
  + \Big( -t_6 \Big) S_{ -2, 1}(x) \\
\nonumber
  & + \Big( t_{11} \Big) S_{ 1, -1}(x)
  + \Big( -t_4 - t_9 +t_{10} +t_{13} \Big) S_{ 0, 0}(x)
  + \Big( -t_8 \Big) S_{ -1, 1}(x) \\
\nonumber
  & + \Big( -t_{11} + t_{12} + t_{15} \Big) S_{ 1, 0}(x)
  + \Big( -t_{10} + t_{14} \Big) S_{ 0, 1}(x)
  + \Big( -t_{12} + t_{16} \Big) S_{ 1, 1}(x) \\
\nonumber
  & \Big( -t_{13} \Big) S_{ -2, -1}(xq)
  + \Big( -t_{15} \Big) S_{ -1, -1}(xq)
  + \Big( -x q t_5 - t_{14} \Big) S_{ -2, 0}(xq) \\
\nonumber
  & + \Big( 0 \Big) S_{ 0, -1}(xq)
  + \Big( -x q t_7 - t_{16} \Big) S_{ -1, 0}(xq)
  + \Big( -x q^2 t_6 \Big) S_{ -2, 1}(xq) \\
\nonumber
  & + \Big( 0 \Big) S_{ 1, -1}(xq)
  + \Big( -x^2 q^2 t_1- x q t_9 \Big) S_{ 0, 0}(xq)
  + \Big( - x q^2 t_8 \Big) S_{ -1, 1}(xq) \\
\nonumber
  & + \Big( -x^2 q^3 t_3- x q t_{11} \Big) S_{ 1, 0}(xq)
  + \Big( -x^2 q^2 t_2- x q^2 t_{10} \Big) S_{ 0, 1}(xq) \\
\label{fneqExampleBig}
  & + \Big( -x^2 q^3 t_4- x q^2 t_{12} \Big) S_{ 1, 1}(xq)
  = 0
\end{align}}
Guided by \eqref{fneq_A_PNAS},
we want the $S_{-2, -1}(x)$, $S_{-2, -1}(xq)$,
$S_{-1, -1}(x)$, $S_{-1, -1}(xq)$, $S_{0, 0}(x)$, and $S_{0, 0}(xq)$ terms to survive,
and the rest of the terms to vanish in \eqref{fneqExampleBig}.
So, we identify the coefficients of unwanted terms with zero,
and solve the resulting linear system for $t_j$'s.
The solution is
\begin{align}
\nonumber
  & t_2 = t_3 = t_5 = t_6 = t_8 = t_9 = t_{10} = t_{11} = t_{14} = 0, \\
\nonumber
  & \qquad t_1, t_7, t_{13} \textrm{ free parameters }, \\
\label{solnExampleBigInterm}
  & t_4 = t_7, \;
  t_{12} = -x q t_7, \;
  t_{13} = t_13, \;
  t_{15} = x q t_7, \;
  t_{16} = -x q t_7.
\end{align}
Applying the solution to \eqref{fneqExampleBig} yields
\begin{align}
\nonumber
  & \Big( t_1 \Big) S_{ -2, -1}(x)
  + \Big( -t_1 + t_7 \Big) S_{ -1, -1}(x)
  + \Big( -t_7 + t_{13} \Big) S_{ 0, 0}(x) \\
\label{solnExampleBigParametric}
  & + \Big( -t_{13} \Big) S_{ -2, -1}(xq)
  + \Big( -x q t_7 \Big) S_{ -1, -1}(xq)
  + \Big( - x^2 q^2 t_{1}\Big) S_{ 0, 0}(xq)
  = 0.
\end{align}
Finally, setting $(t_1, t_7, t_{13})=$ $(1,0,0)$, $(0,1,0)$, $(0,0,1)$,
we recover \eqref{fneq_A_PNAS}.

In general, we start with a single instance of \eqref{srMainDouble}.
We choose and fix intervals
for the first and the second indices \eqref{listIndices1}-\eqref{listIndices2},
and list all instances of \eqref{fneqPrimary} in which
all indices fall in their respective ranges.
Then, we multiply the first equation in that list with $t_1$,
the second one with $t_2$, etc. and add those all up.
This is exactly what we did for our example above.

Unlike \eqref{fneq_A_PNAS}, we do not know how many
other series are needed to uniquely define all of those series at once.
In other words, we do not know if the series we chose to study
is one of the two, or three, or more similar series
all related in that many functional equations.
We need as many functional equations as the number of series
to expect a unique determination.
Therefore, we repeat the following for any combination of the pairs of indices.
We need to keep in mind that all indices but $C_1$ and $C_2$ are preserved in \eqref{fneqPrimary}.
Given any combination of pairs of indices
$A = \Big\{ ( a_1, b_1 ), \ldots, (a_r, b_r) \Big\}$,
we go over the coefficient of each $S_{C_1, C_2}(x)$ and $S_{C_1, C_2}(xq)$
in the counterpart of \eqref{fneqExampleBig}.
If $(C_1, C_2) \in A$, then we discard the coefficient momentarily.
Otherwise, i.e. $(C_1, C_2) \not\in A$,
we equate its coefficient to zero, and append it as another equation
in a linear system in $t_j$'s associated with the index pair subset $A$.
We solve that linear system and plug in the solution in the counterpart of \eqref{fneqExampleBig}.

At this point, we have a reduced version of \eqref{fneqPrimary}
in which all pairs of indices belong to the index pair subset $A$.
The coefficients possibly depend on some of the $t_j$'s.
A setup is possible for a little row reduction,
and we check if we arrive at the following.
{\allowdisplaybreaks
\begin{align}
\nonumber
  & S_{a_1, b_1}(x) = \sum_{ (\alpha, \beta) \in A} f^1_{\alpha, \beta}(x) S_{\alpha, \beta}(x q^\gamma) \\
\nonumber
  & S_{a_2, b_2}(x) = \sum_{ (\alpha, \beta) \in A} f^2_{\alpha, \beta}(x) S_{\alpha, \beta}(x q^\gamma) \\
\nonumber
  & \qquad \vdots \\
\label{fneqGeneral}
  & S_{a_r, b_r}(x) = \sum_{ (\alpha, \beta) \in A} f^r_{\alpha, \beta}(x) S_{\alpha, \beta}(x q^\gamma)
\end{align}}
Equations for some $S_{a_j, b_j}(x)$ may be missing.
This may be because their coefficients in the counterpart of \eqref{fneqExampleBig}
are zero to start with, or become zero after we apply the solution of
the linear system annihilating the unwanted coefficients,
or because the counterpart of \eqref{solnExampleBigParametric}
does not have enough number of parameters,
or it does not give way to an enough number of linearly independent
coefficients for $S_{a_j, b_j}(x)$'s (not for $S_{a_j, b_j}(x q^\gamma)$'s).

You will have noticed that in the heavily guided example above,
the intermediate linear system in $t_j$'s we solved did not have more variables than equations,
so the solution we have purely coincidentally had some free parameters.
We can ensure more variables than equations,
but we cannot ensure a system of functional equations
with a specified number of series in it.
The main obstacle is that we are solving
an auxiliary homogeneous system of linear equations
and applying its solutions in another bunch of expressions.
This is not only requiring some free variables in solutions,
but also asking the free variables to be specified ones.
In fact, we are requiring more than that,
but even this much cannot be guaranteed
with the methods we are employing here.  

We must also emphasize that there are existence results 
for functional equations satisfied by \eqref{srMainDouble}, 
even construction algorithms~\cite{Sp12}.  
However, our motivation is to make a partition theoretic connection.  

In~\cite{Chern}, systems of functional equations for generating functions of 
linked partition ideals~\cite{A_thebluebook} are derived, 
and connections to functional equations satisfied by instances of \eqref{srMainDouble} 
are established.  
Linear combinations of~\eqref{fneqPrimary} are not mentioned.  
In~\cite{KR_20}, instances of \eqref{fneqGeneral} are accepted as input, 
and it is verified that they can be written as linear combinations of \eqref{fneqPrimary}.  

\section{A not so obvious partition generating function}
\label{secFirstApplication}

In this section, 
we will exemplify an interpretation of a double series in the form of \eqref{srMainDouble} 
as a partition generating function.  
The series along with a system of functional equations 
are among the hundreds of outputs of the algorithm outlined in Section \ref{secSrchAlgo}.  
The algorithm was put to a brute force search for parameters in specified intervals.  
At first, the search was adjusted so that 
the series which appear in at least one system of functional equations 
having two, three or four series were determined.  
In other words, if the program found one system, 
it recorded the parameters and moved on.  
Then, for some of those parameters, the algorithm was run one more time 
to list all possible systems of functional equations. 
Those systems were visually inspected for suitability 
for a partition generating function interpretation.  

One would also have liked to have a partition identity in Theorem \ref{thmFirstApplication}, 
but \eqref{eqFirstApplication} does not have a nice 
infinite product representation.  

There are 20 partitions of $n = 14$ satisfying the condition in Theorem \ref{thmFirstApplication}:  
\begin{align*}
  \begin{array}{cc}
  1 & 6 \\ 1 & 6
  \end{array}, \qquad 
  & \begin{array}{cc}
  2 & 5 \\ 2 & 5
  \end{array}, \qquad 
  \begin{array}{ccc}
  1 & & \\ 1 & 3 & 9
  \end{array}, \qquad 
  \begin{array}{ccc}
  1 & & \\ 1 & 4 & 8
  \end{array}, \qquad 
  \begin{array}{ccc}
  1 & & \\ 1 & 5 & 7
  \end{array}, \\[3mm]
  & \begin{array}{ccc}
  2 & & \\ 2 & 4 & 6
  \end{array}, \qquad 
  \begin{array}{ccc}
    & 3 & \\ 1 & 3 & 7
  \end{array}, \qquad 
  \begin{array}{ccc}
    & & 5 \\ 1 & 3 & 5
  \end{array}, 
\end{align*}
along with the 12 Rogers-Ramanujan partitions of $14$.  

\begin{proof}[proof of Theorem \ref{thmFirstApplication}]
 We set the parameters 
 $(B_{11}, B_{22}, B_{12}, D_1, D_2, K_1, K_2, \gamma)$ 
 $=( 6, 2, 2, 2, 1, 2, 1, 1)$, 
 and run the program we described in Section \ref{secSrchAlgo} 
 for $C_1 \in \{ -4, -2, 0, 2, 4 \}$, $C_2 \in \{ -1, 0, 1, 2 \}$, 
 set to list systems of functional equations with four series in them.  
 One of the outputs is the following.  
 {\allowdisplaybreaks
 \begin{align}
 \nonumber 
  S_{-4, -1}(x) & = S_{-4, -1}(xq) + xq S_{-2, 0}(xq) + x^2 q^2 S_{0, 0}(xq) \\
 \nonumber 
  S_{-2, -1}(x) & = S_{-4, -1}(xq) + xq S_{-2, 0}(xq) \\
 \nonumber 
  S_{-2, 0}(x) & = S_{-4, -1}(xq) \\
 \label{eqFirstAppProofFuncEq}
  S_{0, 0}(x) & = S_{-2, -1}(xq) 
 \end{align}}
 It is easy to verify that these functional equations, 
 along with the initial values
 \begin{align}
 \nonumber 
  S_{-4, -1}(0) = S_{-2, -1}(0) = S_{-2, 0}(0) = S_{0, 0}(0) = 1
 \end{align}
 determine the four series uniquely.  
 The initial values correspond to the empty partition of zero.  
 Since the exponent of $x$ keeps track of the number of parts, 
 replacing $x$ by $xq$ adds 1 to all parts in the generated partitions.  
 Multiplying by $xq$ appends a 1 to the partitions, 
 and by $x^2 q^2$ a $1+1$.  
 Then, it is clear that the generated partitions 
 have parts that are repeated at most twice.  
 
 We will argue for all four series at once.  
 We claim that 
 \begin{enumerate}[(i)]
  \item $S_{-4, -1}(x)$ generates partitions 
   with parts which do not repeat more than twice, 
   such that part sizes are at least 2 apart, 
   and repeating part sizes are at least three apart.  
  \item $S_{-2, -1}(x)$ generates partitions listed by 
   $S_{-4, -1}(x)$, but additionally, repeated parts are at least two.  
  \item $S_{-2, 0}(x)$ generates partitions listed by 
   $S_{-4, -1}(x)$, but additionally, no part is equal to 1.  
  \item $S_{0, 0}(x)$ generates partitions listed by 
   $S_{-4, -1}(x)$, but additionally, 
   non-repeated parts are at least two, and repeated parts are at least three.  
 \end{enumerate}
 
 The main argument of the proof is structural induction.  
 We prove the first of \eqref{eqFirstAppProofFuncEq}, 
 the second is shown likewise.  
 The third and the fourth are obvious.  
 The partitions allegedly generated by $S_{-4, -1}(x)$ 
 can be set-partitioned into three classes.  
 
 $\mathcal{U}_0 = $ partitions generated by $S_{-4, -1}(x)$ which have no $1$'s.  
 
 $\mathcal{U}_1 = $ partitions generated by $S_{-4, -1}(x)$ which have a single $1$.  
 
 $\mathcal{U}_{11} = $ partitions generated by $S_{-4, -1}(x)$ which have $1+1$.  
 
 By the inductive hypothesis, $\mathcal{U}_0$ is generated by $S_{-2, 0}(x)$, 
 which is equal to $S_{-4, -1}(xq)$.  
 This is the first term in the first equation of \eqref{eqFirstAppProofFuncEq}.  
 If we delete the 1 from each partition in $\mathcal{U}_1$, 
 by the inductive hypothesis, in the partitions that remain, 
 the parts must be at least three, and three may appear once, or may repeat.  
 These are generated by $S_{-4, -1}(xq^2)$, 
 which in turn is equal to $S_{-2, 0}(xq)$.  
 Along with the $xq$ accounting for the deleted 1, 
 we have the second term $xq S_{-2, 0}(xq)$.  
 For $\mathcal{U}_{11}$, by the inductive hypothesis, 
 the larger parts are at least 3 if they do not repeat, 
 and at least 4 if they do.  
 This is given by $S_{-2, -1}(xq^2)$, which also equals $S_{0, 0}(xq)$.  
 Along with $x^2 q^2$ for the $1+1$, we get the last term $x^2 q^2 S_{0, 0}(xq)$.  
 
 The remaining details are the basis case, 
 which is the empty partition being generated by all four series, 
 and the observation that for non-empty partitions, 
 each step of the structural induction decreases each part, 
 in particular, the largest part by 1.  
\end{proof}

Backed up by Theorem \ref{thmFirstApplication}, 
we can try and give a base partition 
and moves interpretation of \eqref{eqFirstApplication}~\cite{K_19}.  
The index $m$ keeps track of the pairs, 
and $n$ keeps track of the singletons.  
The moves are easy to imagine.  
The singletons are incremented or decremented by one, 
recorded by $(q; q)_n$ in the denominator.  
They are kept at least two apart.  
The pairs move together, such as 
\begin{align}
\nonumber 
  \begin{array}{c}
   \mathbf{a} \\ \mathbf{a} 
  \end{array}
  \longrightarrow
  \begin{array}{c}
   \mathbf{a+1} \\ \mathbf{a+1} 
  \end{array}, 
  \qquad \textrm{ or } \qquad 
  \begin{array}{c}
   \mathbf{a} \\ \mathbf{a} 
  \end{array}
  \longrightarrow
  \begin{array}{c}
   \mathbf{a-1} \\ \mathbf{a-1} 
  \end{array}, 
\end{align}
hence each move alters the weight by 2.  
This explains the $(q^2; q^2)_m$ in the denominator.  
When a pair confronts a singleton, the moves are 
\begin{align}
\nonumber 
  \begin{array}{ccc}
   \mathbf{a} & & \\ \mathbf{a} & & \mathrm{a+2}
  \end{array}
  \longrightarrow
  \begin{array}{ccc}
    & & \mathbf{a+2} \\ \mathrm{a} & & \mathbf{a+2} 
  \end{array}, 
  \qquad \textrm{ or } \qquad 
  \begin{array}{ccc}
    & & \mathbf{a} \\ \mathrm{a-2} & & \mathbf{a}  
  \end{array}
  \longrightarrow
  \begin{array}{ccc}
   \mathbf{a-2} & & \\ \mathbf{a-2} & & \mathrm{a}
  \end{array}.  
\end{align}
For convenience, the singletons or pairs that are being moved 
are shown in boldface.  
We do not allow moves that cause consecutive part sizes, 
or bring the pairs too close together.  
For example; neither can the smaller pair move forward, 
nor can the greater pair move backward, in the following configuration.  
\begin{align}
\nonumber 
  \begin{array}{ccccc}
   \mathrm{a-2} & & & & \mathrm{a+2} \\ 
   \mathrm{a-2} & & \mathrm{a} & & \mathrm{a+2}
  \end{array}
\end{align}
If only the base partition was 
\begin{align}
\label{ptnFirstWrongBase} 
  \begin{array}{ccccccccccccccc}
   1 & & 4 & & & & (3m-2) & & & & & & & & \\
   1 & & 4 & & \cdots & & (3m-2) & & 3m & & (3m+2) & & \cdots & & (3m+2n-2)
  \end{array}, 
\end{align}
the base partition and moves explanation would be almost complete.  
The petty details are easy to supply~\cite{K_19}.  
However; \eqref{ptnFirstWrongBase} has weight 
$6\binom{m+1}{2} + 2\binom{n+1}{2} + 3mn - 4m - 2n$.  
The trained eye would catch that the coefficient of the cross-multiplied term 
is 3, as opposed to the expected 2 in \eqref{eqFirstApplication}.  

This is because the partition \eqref{ptnFirstWrongBase} 
does not have minimal weight among partitions generated by $S_{-4, -1}(x)$.  
For instance, one can check that the partition 
$\begin{array}{cccc} 1 & 4 & & \\ 1 & 4 & 6 & 8 \end{array}$ 
has weight 24 while 
$\begin{array}{cccc} 1 & & 5 & \\ 1 & 3 & 5 & 7 \end{array}$ 
has weight 22.  The most convenient way to explain this is 
to allow the smallest singleton $3m$ 
in the na\"{\i}ve guess \eqref{ptnFirstWrongBase} of a base partition 
to make some more backward moves.  
\begin{align}
\nonumber 
 \begin{array}{c} 1 \\ 1 \end{array}
 \begin{array}{c} 4 \\ 4 \end{array}
 \begin{array}{c}  \\ \cdots \end{array}
 \begin{array}{c} \mathrm{(3m-2)} \\ \mathrm{(3m-2)} \end{array}
 \begin{array}{c}  \\ \mathbf{3m} \end{array}
 \begin{array}{c}  \\ \mathrm{(3m+2)} \end{array}
 \begin{array}{c}  \\ \cdots \end{array}
 \begin{array}{c}  \\ \mathrm{(3m+2n-2)} \end{array}
\end{align}
\begin{align}
\nonumber 
  \Big\downarrow \textrm{ one backward move }
\end{align}
\begin{align}
\nonumber 
 \begin{array}{c} 1 \\ 1 \end{array}
 \begin{array}{c} 4 \\ 4 \end{array}
 \begin{array}{c}  \\ \cdots \end{array}
 \underbrace{\begin{array}{c} \mathrm{(3m-2)} \\ \mathrm{(3m-2)} \end{array}
 \begin{array}{c}  \\ \mathbf{(3m-1)} \end{array}}_{\mathrm{ ! }}
 \begin{array}{c}  \\ \mathrm{(3m+2)} \end{array}
 \begin{array}{c}  \\ \cdots \end{array}
 \begin{array}{c}  \\ \mathrm{(3m+2n-2)} \end{array}
\end{align}
\begin{align}
\nonumber 
  \Big\downarrow \textrm{ adjustment }
\end{align}
\begin{align}
\nonumber 
 \begin{array}{c} 1 \\ 1 \end{array}
 \begin{array}{c} 4 \\ 4 \end{array}
 \begin{array}{c}  \\ \cdots \end{array}
 \begin{array}{c} \mathrm{(3m-5)} \\ \mathrm{(3m-5)} \end{array}
 \begin{array}{c}  \\ \mathbf{(3m-3)} \end{array}
 \begin{array}{c} \mathrm{(3m-1)} \\ \mathrm{(3m-1)} \end{array}
 \begin{array}{c}  \\ \mathrm{(3m+2)} \end{array}
 \begin{array}{c}  \\ \cdots \end{array}
 \begin{array}{c}  \\ \mathrm{(3m+2n-2)} \end{array}
\end{align}
The backward move decreases the weight by one, 
and the adjustment preserves it.  
Afterwards, we see that a backward move for each of the 
greater singletons is possible, 
yielding the intermediate partition
\begin{align}
\nonumber 
 \begin{array}{c} 1 \\ 1 \end{array}
 \begin{array}{c} 4 \\ 4 \end{array}
 \begin{array}{c}  \\ \cdots \end{array}
 \begin{array}{c} \mathrm{(3m-5)} \\ \mathrm{(3m-5)} \end{array}
 \begin{array}{c}  \\ \mathbf{(3m-3)} \end{array}
 \begin{array}{c} \mathrm{(3m-1)} \\ \mathrm{(3m-1)} \end{array}
 \begin{array}{c}  \\ \mathrm{(3m+1)} \end{array}
 \begin{array}{c}  \\ \cdots \end{array}
 \begin{array}{c}  \\ \mathrm{(3m+2n-3)} \end{array}.  
\end{align}
Together with the initial backward move on the smallest singleton, 
this decreases the weight by $m$.  
This maneuver is possible over all but the smallest pair, 
decreasing the weight by $(m-1)n$.  
At that point, no more backward moves are possible on any 
pair or singleton without losing parts 
or without having negative parts.  
This constructs the base partition, 
and its weight is $6\binom{m+1}{2} + 2\binom{n+1}{2} + 2mn - 4m - n$, 
as desired.  
The initial movement of singletons in the forward phase 
may involve adjustments, but they use the same idea as above.  
We write only the crucial steps for the sake of brevity.  

Variants of the above construction are possible.  
For instance, the series 
\begin{align}
\nonumber 
\sum_{m, n \geq 0} f(m, n) x^m q^n 
= \sum_{m, n} \frac{ q^{ 6\binom{m+1}{2} + 2\binom{n+1}{2} + 2mn} x^{2m+n} }
  { ( q^2; q^2)_m ( q; q)_n }
\end{align}
generates partitions into parts that repeat at most twice, 
such that singletons  are at least two apart, 
part sizes of pairs are at least three apart, 
singletons are at least two, 
part sizes of pairs are at least three, 
a singleton may follow a pair, 
but a pair must be preceded by two unused part sizes.  
For instance, 
\begin{align}
\nonumber 
  \begin{array}{cc}
   3 & 6 \\ 3 & 6
  \end{array}, \qquad 
  \begin{array}{cccc}
   3 & & 7 & \\ 3 & 4 & 7 & 8
  \end{array}, \qquad 
  \begin{array}{cc}
    & 5 \\ 2 & 5
  \end{array}
\end{align}
are among the generated partitions.

In conclusion, the partition class in Theorem \ref{thmFirstApplication} was described 
with the help of the selected system of functional equations.  
Once the partition class was known, 
the series turned out to be fit for explanation 
via a base partition and moves framework~\cite{K_19}.  
We must impress that the base partition and the moves framework 
are not obvious unless one knows what partitions are generated.  
The moves could be guessed, but the base partition is tricky to construct.

\section{Narrowing the search space}
\label{secNarrowingSrchSpc}

The careful reader will have noticed that in Section \ref{secFirstApplication}, 
we said $C_1 \in \{ -4, -2, 0, 2, 4 \}$ 
rather than $ -4 \leq C_1 \leq 4 $.  
The reason for this is the pairwise differences between 
the first indices listed in \eqref{listIndices1} are 
\begin{align}
\nonumber 
  K_1, \quad B_{11}-\gamma D_1, \quad B_{12} - \gamma D_1, \quad \gamma D_1.  
\end{align}
Thus, once the parameters $B_{11}$, $B_{22}$, $B_{12}$, 
$D_1$, $D_2$, $K_1$, $K_2$, and $\gamma$ are chosen and fixed, 
instances of functional equations listed in \eqref{fneqPrimary} 
cannot relate instances of series \eqref{srMainDouble} 
unless the first indices differ by a multiple of 
\begin{align}
\nonumber 
  d_1 = \mathrm{g.c.d.}\left(K_1, B_{11}, B_{12}, \gamma D_1 \right)
\end{align}
Likewise, functional equations listed in \eqref{fneqPrimary} 
can only relate instances of series \eqref{srMainDouble} 
when the second indices differ by a multiple of 
\begin{align}
\nonumber 
  d_2 = \mathrm{g.c.d.}\left(K_2, B_{22}, B_{12}, \gamma D_2 \right)
\end{align}
One can verify that $d_1 = 2$ and $d_2 = 1$ 
in Theorem \ref{thmFirstApplication}.  

This observation allows us to give a more accurate lower bound on 
the number of series that appear, and that are actually related, 
by the functional equations in \eqref{fneqPrimary}.  
In particular, 
\begin{align}
\label{eqIndexPairCount} 
  \textrm{ the number of instances of \eqref{fneqPrimary} }
  \leq 2 \left( \frac{ \Delta m_1 }{d_1} + 1 \right) 
   \left( \frac{ \Delta m_2 }{d_2} + 1 \right),  
\end{align}
where $\Delta m_j = M_j - m_j$, 
and $\left[ m_j, M_j \right]$ are the intervals for $C_j$, $j = 1, 2$.  
It is not strictly necessary, but convenient, to make the assumption 
that $d_j \vert \Delta m_j$.  
The doubling is required because we take 
series with argument $x$ and $x q^{\gamma}$ independent.  
The inequality instead of equation in \eqref{eqIndexPairCount} 
is exemplified in \eqref{fneqExampleBig}, 
where some coefficients are zero.  
That means, the corresponding series do not appear in 
any of the functional equations.  
It is also possible to quantify and incorporate the error so that
\eqref{eqIndexPairCount} is made an equality 
for $\Delta m_j$ greater than effectively calculable thresholds.

On the other hand, for each of the equation types in \eqref{fneqPrimary}, 
we will calculate the maximum difference among the 
appearing first indices, and the same for the second indices.  
Call those maximum differences in the first indices $x_j$ and 
in the second indices $y_j$ for $j = 1, 2, 3$.  
It is straightforward to see that 
\begin{align}
\nonumber 
  x_1 = \mathrm{max}\left\{ K_1, 
    \left\vert B_{11} - \gamma D_1 \right\vert, 
    \left\vert K_1 - B_{11} + \gamma D_1 \right\vert \right\}, 
\qquad
  y_1 = \left\vert B_{12} - \gamma D_1 \right\vert, 
\end{align}
\begin{align}
\nonumber 
  x_2 = \left\vert B_{12} - \gamma D_2 \right\vert, 
\qquad 
  y_2 = \mathrm{max}\left\{ K_2, 
    \left\vert B_{22} - \gamma D_2 \right\vert, 
    \left\vert K_1 - B_{22} + \gamma D_2 \right\vert \right\},
\end{align}
\begin{align}
\nonumber 
  x_3 = \gamma D_1, 
\qquad 
  y_3 = \gamma D_2.  
\end{align}
This will enable us to find the number of instances of 
each type of equation in \eqref{fneqPrimary}.  
Visually, observe that the indices in the first equation of  \eqref{fneqPrimary}
are lattice points on $\mathbb{Z}^2$:
\begin{align}
\nonumber
  \left\{ ( C_1,  C_2 ), (C_1 + K_1, C_2),
   (C_1 + B_{11} - \gamma D_1, C_2 + B_{12} - \gamma D_2) \right\}.
\end{align}
The tightest rectangle containing these points
has width $x_1$ and height $y_1$ as given above.
\begin{center}
\includegraphics[scale=0.1]{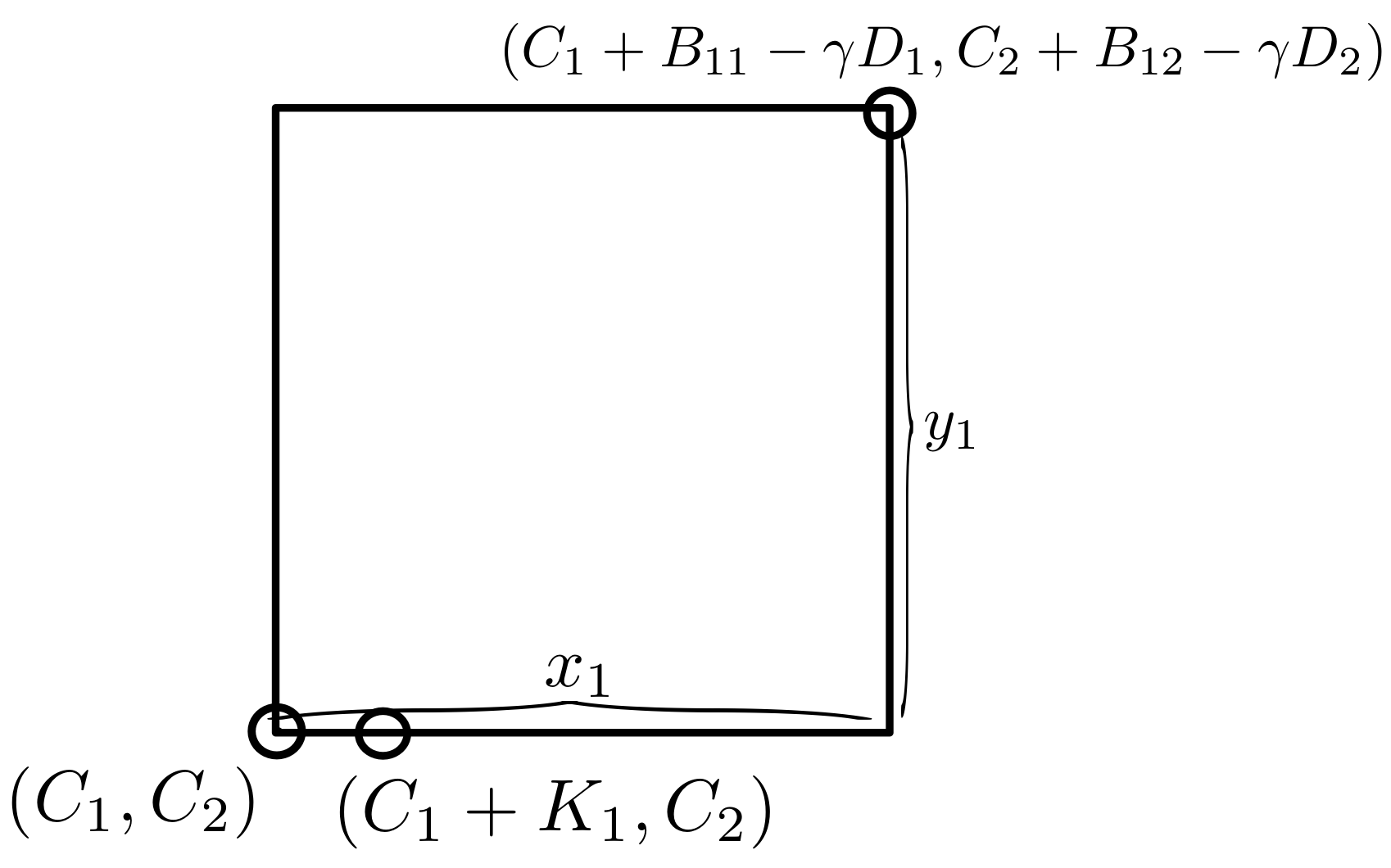}
\end{center}
Similar pictures can be drawn for $x_2$ and $y_2$, or $x_3$ and $y_3$.
It is possible that the rectangles are degenerate,
i.e. the height or the width is zero.

In brief,  we will be able to calculate the number of equations
in the counterparts of the list \eqref{fneqRunningExample}.  
We will have a total of 
\begin{align}
\label{eqFneqCount}
 \sum_{j = 1}^3 \left( \frac{ \Delta m_1 - x_j}{d_1} + 1 \right) 
  \left( \frac{ \Delta m_2 - y_j}{d_2} + 1 \right)
\end{align}
equations.  
Therefore, we will have a total of that many variables ($t$'s) in 
counterparts of \eqref{fneqExampleBig}.  

If we are searching for systems of functional equations 
relating exactly $d$ instances of \eqref{srMainDouble}, 
then we will keep $2d$ of the series 
the number of which is bounded by \eqref{eqIndexPairCount} 
on the side, equate the remaining to zero, 
and solve for the $t$'s.  
These equations are guaranteed to have non-trivial solutions when 
\begin{align}
\label{ineqMainDeltaM}
 \sum_{j = 1}^3 \left( \frac{ \Delta m_1 - x_j}{d_1} + 1 \right) 
  \left( \frac{ \Delta m_2 - y_j}{d_2} + 1 \right)
 > 2 \left( \frac{ \Delta m_1 }{d_1} + 1 \right) 
   \left( \frac{ \Delta m_2 }{d_2} + 1 \right) - 2d.  
\end{align}
This is possible to ensure as long as we choose $\Delta m_1$ and $\Delta m_2$ 
proportional and large enough.  
In that case, both sides of \eqref{ineqMainDeltaM} become 
quadratic expressions, and the leading term on the left hand side 
is $\frac{3}{2}$ times the leading term on the right hand side.  

However, one cannot guarantee a solution in the desired form 
when we apply the non-trivial solution to the $2d$
coefficients we kept on the side.  
After all, we are applying the solution of one homogenous linear system 
to another set of expressions.  
This is, again, why the procedure is a search algorithm at best.  

A brute force search on parameter lists is infeasible, 
as pointed out in Section \ref{secFirstApplication}.  
But, the numbers on either side of \eqref{ineqMainDeltaM} 
can be calculated much faster than constructing and solving a linear system.  
We also chose to discard instances of \eqref{srMainDouble} 
in which $\mathrm{g.c.d.}(B_{11}, B_{22}, B_{12}, K_1, K_2) > 1$ 
because they are dilations of series with strictly smaller parameters.  
One can of course incorporate $C_j$'s to this picture, 
and ensure that the series is not a dilation.  
But, we did not go into that fine tuning.  

We restricted our attention to systems in which
both counts in \eqref{ineqMainDeltaM} are at most 200, 
and the parameters $B_{11}$, $B_{22}$, $B_{12}$, $D_1$, $D_2$, $K_1$, $K_2$ and $\gamma$
are within small bounds.  
In particular, we restricted 
$1 \leq B_{\cdots} \leq 8$ and $1 \leq D_\cdot, K_\cdot, \gamma \leq 4$.  
We had hundreds of parameter lists to look at.
Given that the list of systems of functional equations
still needs non-automated checking,
this is still too many.
Besides, one would like to have a partition identity,
or a potential partition identity so that the partition generating function is more appealing.

At this point, for each of the hundreds of parameter lists,
we searched for periodic infinite product representations
with the help of Euler's algorithm~\cite{A_anotherbook, KR_15}.
For each fixed set of parameters in \eqref{srMainDouble},
we substituted $x = 1$, found Maclaurin coefficients
\begin{align}
\nonumber
  S_{C_1, C_2}(1) = 1 + b_1 q + b_2 q^2 + \cdots + b_M q^M + O(q^{M+1}),
\end{align}
for $(C_1, C_2)$ $ \in \left[ -B_{11}, \leq B_{11}\right] \times \left[ -B_{22}, \leq B_{22}\right]$, 
found the first $M$ exponents in the infinite product representation
\begin{align}
\nonumber
  S_{C_1, C_2}(1) = \frac{ (1 + O(q^{M+1})) }{ (1 - q)^{a_1} (1 - q^2)^{a_2} \cdots (1 - q^M)^{a_M} },
\end{align}
and checked if there is any periodicity with period $k$.
We used $M=50$ and $k \leq 24$. 
We discarded the series for which $S_{C_1, C_2}(0) \neq 1$.  

This brings an acceptable redundancy as far as the parameters $D_1$ and $D_2$ are concerned,
but the redundancy is well worth the results.
Coupled with this series-product identity search,
we now have less than a hundred series to look at.
Plus, for each fixed list of parameters
$(B_{11}, B_{22}, B_{12}, D_1, D_2, K_1, K_2, \gamma)$,
the number of pairs $(C_1, C_2)$ gave us an idea of
how many functional equations there will be
in a system of functional equations that uniquely define the $S_{\cdot, \cdot}(x)$ at hand.
Along with some classical and more recent results
(Andrews-Gordon identities for $k = 3$~\cite{A_PNAS}, 
Kanade-Russell mod 9 conjectures~\cite{KR_15, K_19_KR}, 
many apparent cases of Euler's formula that's a corollary of 
the $q$-binomial theorem~\cite{A_thebluebook})
we found Theorems \ref{thmYalcinGenFunc} and \ref{thmYalcinGenFunc2},
and gave combinatorial explanations.

The series in consideration has the fixed parameters 
$( B_{11}, B_{22}, B_{12}, D_1, D_2, K_1, K_2, \gamma )$ 
$=( 2, 1, 1, 2, 1, 1, 1, 1 )$.  

\begin{theorem}
\label{thmYalcinGenFunc}
 Let $t_1(m, n)$ be the number of partitions of $n$ 
 into $m$ bicolored (red and blue) parts 
 satisfying the following conditions.  
 \begin{enumerate}
  \item Neither blue nor red parts repeat.   
  \item when we go through parts from the smallest to largest, 
   for each red part $b$, there must be a blue $b$ or a blue $b+1$. 
   Moreover, these blue parts should be distinct for each $\textcolor{red}{b}$.  
 \end{enumerate}
 Let $t_2(m, n)$ be the number of partitions 
 enumerated by $t_1(m, n)$ in which 
 a red 1 and a blue 1 do not appear together.  
 Then, 
 \begin{align}
 \nonumber 
  \sum_{m, n \geq 0} t_1(m, n) x^m q^n = S_{0,0}(x) 
  & = \sum_{m, n \geq 0} \frac{ q^{ 2 \binom{m+1}{2} 
      + \binom{n+1}{2} + mn } x^{2m+n} }
    { (q; q)_m (q; q)_n },  \\
 \nonumber 
  \sum_{m, n \geq 0} t_2(m, n) x^m q^n = S_{1,0}(x) 
  & = \sum_{m, n \geq 0} \frac{ q^{ 2 \binom{m+1}{2} 
      + \binom{n+1}{2} + mn + m } x^{2m+n} }
    { (q; q)_m (q; q)_n }.  
 \end{align}
\end{theorem} 

The four partitions of $n = 4$ satisfying the conditions set forth by $t_1(m, n)$ are 
\begin{align}
\nonumber 
  \textcolor{blue}{4}, \qquad 
  \textcolor{blue}{1} + \textcolor{blue}{3}, \qquad 
  \textcolor{red}{2} + \textcolor{blue}{2}, \qquad 
  \textcolor{red}{1} + \textcolor{blue}{1} + \textcolor{blue}{2}.  
\end{align}
Moreover, $\textcolor{red}{5} + \textcolor{blue}{5}$ is enumerated by $t_1(2, 10)$, 
$\textcolor{red}{1} + \textcolor{blue}{2} + \textcolor{blue}{5} + \textcolor{blue}{14}$ 
by $t_1(4, 22)$; 
while $\textcolor{red}{4} + \textcolor{red}{5} + \textcolor{blue}{5}$ 
is not enumerated by $t_1(3, 14)$, 
and $\textcolor{blue}{1} + \textcolor{red}{6} + \textcolor{blue}{6} + \textcolor{red}{9}$ 
not by $t_1(4, 22)$.  

Before we present the proof, 
let us apply two of Euler's identities~\cite{A_thebluebook} to $S_{0,0}(1)$.  
{\allowdisplaybreaks
\begin{align}
\nonumber 
  & \sum_{m, n \geq 0} \frac{ q^{ 2 \binom{m+1}{2} + \binom{n+1}{2} + mn } }
    { (q; q)_m (q; q)_n }
  = \sum_{ m \geq 0 } \frac{ q^{ 2 \binom{m+1}{2} } }{ (q; q)_m } \; 
   \sum_{n \geq 0} \frac{ q^{ \binom{n+1}{2} } \left( q^m \right)^n }{ (q; q)_n }
  = \sum_{ m \geq 0 } \frac{ q^{ 2 \binom{m+1}{2} } }{ (q; q)_m } \; 
    \left( -q^{m+1}; q \right)_\infty \\ 
\label{eqYalcin3Regular}
  & ( -q; q)_\infty
    \sum_{ m \geq 0 } \frac{ q^{ 2 \binom{m+1}{2} } }{ (q^2; q^2)_m }
  = ( -q; q)_\infty ( -q^2; q^2)_\infty 
  = \frac{ (q^4; q^4)_\infty }{ (q; q)_\infty }
\end{align}}
In other words, $\sum_{m \geq 0} t_1(m, n)$ is also the number of partitions 
of $n$ into parts that repeat at most three times.  
However, $t_1(m, n)$ is not the partitions of $n$ into $m$ parts 
none of which repeat at most three times.  

\begin{proof}[proof of Theorem \ref{thmYalcinGenFunc}]
 One output of the algorithm in Section \ref{secSrchAlgo} 
 is the following system of functional equations.  
 \begin{align}
 \nonumber 
  S_{0,0}(x) & = (1 + x^2 q^2) S_{0,0}(xq) + (xq + x^2 q^3) S_{1, 0}(xq) \\
 \nonumber 
  S_{1,0}(x) & = S_{0,0}(xq) + (xq + x^2 q^3) S_{1, 0}(xq) 
 \end{align}
 Given the initial values $S_{0,0}(0)=S_{1,0}(0)=1$, 
 it is routine to check that the functional equations 
 uniquely determine $S_{0,0}(x)$ and $S_{1,0}(x)$.  
 It is also possible to translate these
 to recurrences between enumerants $t_1(m, n)$ and $t_2(m, n)$~\cite{Yalcin_thesis}.  
 
 Let $\mathcal{T}_1$ be the collection of partitions generated by $S_{0,0}(x)$, 
 and $\mathcal{T}_2$ be that by $S_{1,0}(x)$.  
 Recall that replacing $x$ by $xq$ amounts to adding 1 to each part 
 in all generated partitions in a partition generating function, 
 and multiplying a series by, say, $x^2 q^2$ can be interpreted 
 as appending a $1+1$ at the beginning of each generated partition, 
 we can make the following construction.  
 
 The first of the functional equations written above implies that 
 if $(\lambda_1, \lambda_2, \ldots, \lambda_m)$ $\in \mathcal{T}_1$, 
 then $(\lambda_1+1, \lambda_2+1, \ldots, \lambda_m+1)$ and 
 $(1,1,\lambda_1+1, \lambda_2+1, \ldots, \lambda_m+1)$ $\in \mathcal{T}_1$, too.  
 Also, if $(\lambda_1, \lambda_2, \ldots, \lambda_m)$ $\in \mathcal{T}_2$, 
 then $(1,\lambda_1+1, \lambda_2+1, \ldots, \lambda_m+1)$ and 
 $(1,2,\lambda_1+1, \lambda_2+1, \ldots, \lambda_m+1)$ $\in \mathcal{T}_1$.  
 The second of the functional equations is interpreted similarly.  
 
 Observe that we have not painted the parts yet.  
 The interpretations alluded to in the previous paragraph 
 can be used recursively to construct all partitions 
 in $\mathcal{T}_1$ and $\mathcal{T}_2$.  
 
\begin{tikzpicture}[
  every node/.style={font=\small},
  level 1/.style={sibling distance=25mm},
  level 2/.style={sibling distance=15mm},
  edge from parent/.style={draw,-{Latex}},
  thick
]

\colorlet{blueText}{blue}
\colorlet{redText}{red}
\colorlet{greenBox}{green}

\node[] (root1) {$\emptyset$}
child{node[blueText]{1+\textcolor{redText}{1}}
  child {node[blueText] {2+ \textcolor{redText}{2}}}
  child {node {\textcolor{blueText}{2}+\textcolor{redText}{2}+\textcolor{blueText}{1}+\textcolor{redText}{1}}}
  child {node[rectangle,draw,blueText] {2+ \textcolor{redText}{2}}}};

\node[circle,draw,right=7cm of root1] (root2) {$\emptyset$}
  child {node[blueText] {1}}
  child {node[blueText] {2+ \textcolor{redText}{1}}}
  child {node[rectangle,draw,blueText] {1}}
  child {node[rectangle,draw,blueText] {2+ \textcolor{redText}{1}}};

\end{tikzpicture}
 
 In the above figure, we showed the construction tree with height 1.  
 The roots contain the empty partition $\varepsilon$ of zero.  
 The one on the left is for $\mathcal{T}_1$, 
 and the right one is for $\mathcal{T}_2$
 The enclosed partitions come from $\mathcal{T}_2$.  
 We still need to explain why colors are needed.  
 If we go deeper in these trees, 
 we encounter multiple occurrences of the same partition.  
 To distinguish these multiple generations, 
 we switch to colored partitions.  

\begin{center}
\includegraphics[scale=1.2]{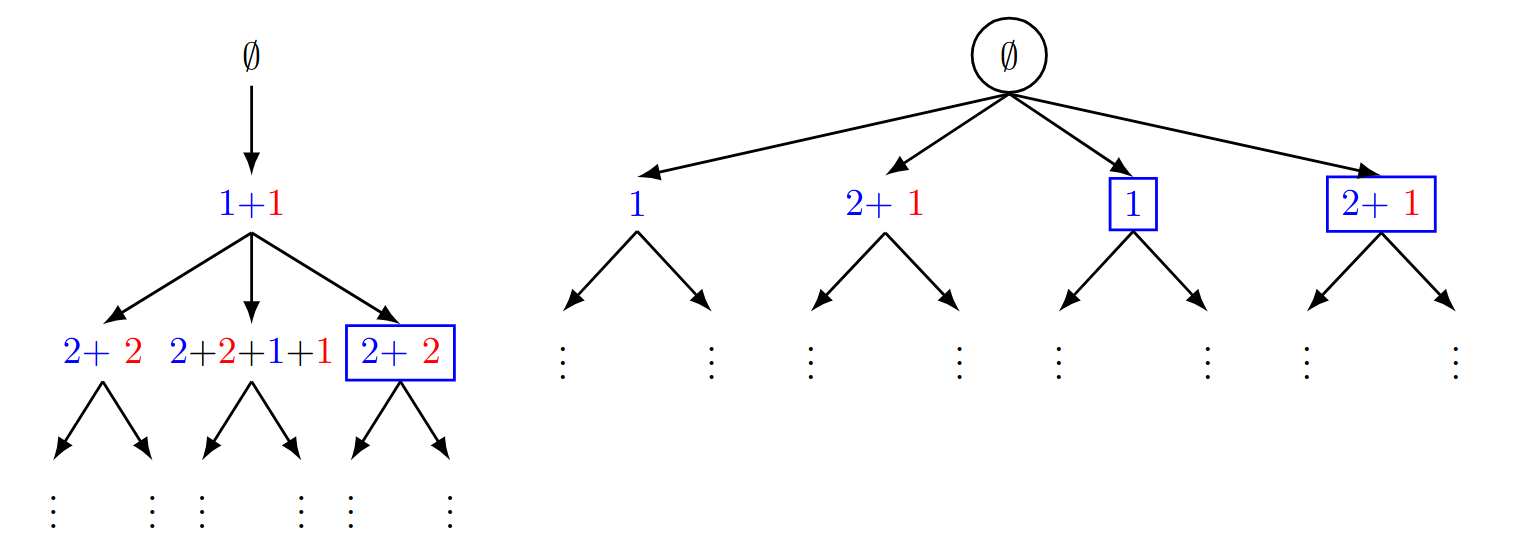}
\end{center}

 After the preparation, the proof is by structural induction.  
 The basis cases are $\varepsilon \in \mathcal{T}_1 \cap \mathcal{T}_2$ 
 and $1 \in \mathcal{T}_1 \cap \mathcal{T}_2$, as seen from the 
 construction tree.  
 
 Suppose the claims of the theorem hold for all partitions 
 in $\mathcal{T}_1$ and $\mathcal{T}_2$ 
 with largest parts $\lambda_m$.  
 We will show that they hold for those 
 with largest part $\lambda_m + 1$, as well.  
 We have the following alternative steps in the construction.  
 $(\lambda_1, \lambda_2, \ldots, \lambda_m)$ $\in \mathcal{T}_1$ implies that 
 \begin{itemize}
  \item[A1.] $(\lambda_1+1, \lambda_2+1, \ldots+1, \lambda_m+1)$ $\in \mathcal{T}_1$, 
  \item[A2.] $(\textcolor{red}{1}, \textcolor{blue}{1}, 
   \lambda_1+1, \lambda_2+1, \ldots+1, \lambda_m+1)$ $\in \mathcal{T}_1$, 
  \item[A3.] $(\lambda_1+1, \lambda_2+1, \ldots+1, \lambda_m+1)$ $\in \mathcal{T}_2$.  
 \end{itemize}
 Similarly, if 
 $(\lambda_1, \lambda_2, \ldots, \lambda_m)$ $\in \mathcal{T}_2$, then 
 \begin{itemize}
  \item[B1.] $(\textcolor{blue}{1}, 
    \lambda_1+1, \lambda_2+1, \ldots+1, \lambda_m+1)$ $\in \mathcal{T}_1$, 
  \item[B2.] $(\textcolor{gray}{1, 2}, 
    \lambda_1+1, \lambda_2+1, \ldots+1, \lambda_m+1)$ $\in \mathcal{T}_1$, 
  \item[B3.] $(\textcolor{blue}{1}, 
    \lambda_1+1, \lambda_2+1, \ldots+1, \lambda_m+1)$ $\in \mathcal{T}_2$, 
  \item[B4.] $(\textcolor{gray}{1, 2}, 
    \lambda_1+1, \lambda_2+1, \ldots+1, \lambda_m+1)$ $\in \mathcal{T}_2$.  
 \end{itemize}
 We used gray to indicate that the colors may be different in different cases.  
 
 First of all, we need to show that none of the steps violate 
 the conditions set forth in the theorem.  
 A1, A2, A3, B1 and B3 clearly do not violate any conditions.  
 For B2 and B4 we need to be more careful, 
 since we cannot have two copies of $\textcolor{red}{2}$ 
 or two copies of $\textcolor{blue}{2}$.  
 We have three cases to consider: 
 \begin{enumerate}[1.]
  \item $\lambda_1 \geq 2$ (regardless of its color), then B2 and B4 take the form 
    $(\textcolor{red}{1}, \textcolor{blue}{2}, 
      \lambda_1+1, \lambda_2+1, \ldots+1, \lambda_m+1)$ $\in \mathcal{T}_i$, 
      $i =$1 or 2 as appropriate.  
  \item $\lambda_1 = \textcolor{blue}{1}$, then B2 and B4 would be applied as 
    $(\textcolor{blue}{1}, \textcolor{red}{2}, 
      \lambda_1+1, \lambda_2+1, \ldots+1, \lambda_m+1)$ $\in \mathcal{T}_i$, 
      $i =$1 or 2 as appropriate.  
      We cannot have another $\textcolor{red}{2}$, 
      since the smallest part $\lambda_1$ is $\textcolor{blue}{1}$.  
  \item $\lambda_1 = \textcolor{red}{1}$, then B2 and B4 would be applied as 
    $(\textcolor{red}{1}, \textcolor{blue}{2}, 
      \lambda_1+1, \lambda_2+1, \ldots+1, \lambda_m+1)$ $\in \mathcal{T}_i$, 
      $i =$1 or 2 as appropriate.  
      This is not a problem, since $\textcolor{blue}{1}$ and $\textcolor{red}{1}$ 
      cannot exits in $\lambda$ together due to the conditions 
      stipulated by $\mathcal{T}_2$.  
 \end{enumerate}
 As a result, all of these steps preserve $\mathcal{T}_1$ and $\mathcal{T}_2$.  
 
 Now we know that this construction, 
 beginning with the empty partition $\varepsilon$ of zero, 
 creates subsets of $\mathcal{T}_1$ and $\mathcal{T}_2$.  
 To conclude the proof, 
 it remains to show that this construction creates 
 all partitions in $\mathcal{T}_1$ and $\mathcal{T}_2$ uniquely.  
 
 Let $(\lambda_1, \ldots, \lambda_{m-1}, \lambda_m + 1)$ 
 be a partition in $\mathcal{T}_i$, $i=$ 1 or 2.  
 We have different cases to consider: 
 \begin{enumerate}[1.]
  \item If $\lambda$ does not contain any $1$'s, 
   then $\lambda$ is obtained from A1 applied 
   to the partition $(\lambda_1-1, \ldots, \lambda_{m-1}-1, \lambda_m)$ 
   $\in \mathcal{T}_1$ for $\lambda \in \mathcal{T}_1$, 
   or A3 to the partition $(\lambda_1-1, \ldots, \lambda_{m-1}-1, \lambda_m)$ 
   $\in \mathcal{T}_2$ for $\lambda \in \mathcal{T}_2$.  
   These are the unique steps, 
   because all others introduce 1 (red or blue) as a part.  
  \item If $\lambda$ contains a $\textcolor{blue}{1}$, 
   but not a $\textcolor{red}{1}$, then we have two cases to consider: 
   \begin{enumerate}[(a)]
    \item If $\lambda$ does not contain $\textcolor{red}{2}$, 
     then it is ontained from B1 applied to the partition 
     $(\lambda_2-1, \ldots, \lambda_{m-1}-1, \lambda_m)$ 
     $\in \mathcal{T}_2$ if $\lambda \in \mathcal{T}_1$, 
     and from B3 applied to the partition 
     $(\lambda_2-1, \ldots, \lambda_{m-1}-1, \lambda_m)$
     $\in \mathcal{T}_2$ if $\lambda \in \mathcal{T}_2$.  
     These are unique ways.  
     We cannot use A1 or A3, since they do not give a 1.  
     We cannot use A2, since it creates a $\textcolor{red}{1}$.  
     We cannot use B2 or B4, either, 
     since they give $\textcolor{red}{1}$ or $\textcolor{red}{2}$.  
    \item If $\lambda$ contains a $\textcolor{red}{2}$, 
     then it must contain a $\textcolor{blue}{2}$ or a $\textcolor{red}{2}$.  
     \begin{enumerate}[i.]
      \item If $\lambda$ contains a $\textcolor{blue}{2}$, 
       then we cannot have used B1 or B3 to obtain $\lambda$.  
       These require applying B1 or B3 to a partition in $\mathcal{T}_2$ 
       which contains both $\textcolor{blue}{1}$ and $\textcolor{red}{1}$, 
       which is impossible.  
       Thus, $\lambda$ is obtained via step B2 applied to 
       $(\lambda_3-1, \ldots, \lambda_{m-1}-1, \lambda_m)$ 
       $\in \mathcal{T}_1$ where $\lambda_3-1 = 1$, 
       if $\lambda \in \mathcal{T}_1$.  
       Similarly, $\lambda$ is obtained via step B4 applied to 
       $(\lambda_3-1, \ldots, \lambda_{m-1}-1, \lambda_m)$ 
       $\in \mathcal{T}_1$ where $\lambda_3-1 = 1$, 
       if $\lambda \in \mathcal{T}_2$.  
      \item If $\lambda$ does not contain $\textcolor{blue}{2}$, 
       then $\lambda$ is obtained via B1 applied to 
       $(\lambda_2-1, \ldots, \lambda_{m-1}-1, \lambda_m)$ 
       where $\lambda_3 - 1 = \textcolor{blue}{2}$ 
       and $\lambda_2 - 1 = \textcolor{red}{1}$, 
       if $\lambda \in \mathcal{T}_1$.  
       Similarly, if $\lambda \in \mathcal{T}_2$, 
       we must have obtained $\lambda$ via applying B3 to 
       $(\lambda_2-1, \ldots, \lambda_{m-1}-1, \lambda_m)$ 
       where $\lambda_3 - 1 = \textcolor{blue}{2}$ 
       and $\lambda_2 - 1 = \textcolor{red}{1}$.  
     \end{enumerate}
   \end{enumerate}
  \item If $\lambda$ contains a $\textcolor{blue}{1}$ 
   and a $\textcolor{red}{1}$, 
   then $\lambda$ must be in $\mathcal{T}_1$ and not in $\mathcal{T}_2$.  
   There is only one way to have obtained $\lambda$, 
   namely using A2 on $(\lambda_3-1, \ldots, \lambda_{m-1}-1, \lambda_m)$.  
   Uniqueness follows from the fact that there is only one step 
   which introduces both $\textcolor{blue}{1}$ and $\textcolor{red}{1}$.  
  \item $\lambda$ contains a $\textcolor{red}{1}$, 
   but not a $\textcolor{blue}{1}$.  
   Then, we must have $\textcolor{blue}{2}$, as well.  
   Thus, if $\lambda \in \mathcal{T}_1$ we must have used 
   B2 on $(\lambda_3-1, \ldots, \lambda_{m-1}-1, \lambda_m)$ 
   where $\lambda_3 - 1 \neq 1$; 
   and if $\lambda \in \mathcal{T}_2$, we must have used 
   B4 on $(\lambda_3-1, \ldots, \lambda_{m-1}-1, \lambda_m)$ 
   where $\lambda_3 - 1 \neq 1$, to obtain $\lambda \in \mathcal{T}_2$.  
   Uniqueness is obvious, since there is only one step 
   which gives $\textcolor{red}{1}$ without $\textcolor{blue}{1}$.  
 \end{enumerate}
 This concludes the proof.  
\end{proof} 

As a result, we start with an otherwise recognizable double series, 
and via a functional equation system discover another interpretation 
which is not clear from the series.  

It is possible to interpret $S_{0,0}(x)$ in Theorem \ref{thmYalcinGenFunc} 
in the (colored) moves framework.  
The base partition is 
\begin{align}
\nonumber 
  \begin{array}{cccccccc}
   \textcolor{blue}{1} & \textcolor{blue}{2} & & \textcolor{blue}{m} & & & & \\ 
   \textcolor{red}{1} & \textcolor{red}{2} & \cdots & 
    \textcolor{red}{m} & (\textcolor{blue}{m+1}) & (\textcolor{blue}{m+1}) & \cdots & 
    (\textcolor{blue}{m+n}) 
  \end{array}.  
\end{align}
The types of parts are singletons, $\textcolor{blue}{b}$, 
which are not paired with any red part, 
and pairs $\begin{array}{cc} & \textcolor{blue}{b} \\ \textcolor{red}{b-1} & \end{array}$
or $\begin{array}{c} \textcolor{blue}{b} \\ \textcolor{red}{b} \end{array}$.  

We describe the backward moves on pairs first.  
We have two cases to consider.  
\begin{align}
\nonumber 
  \begin{array}{cc} & \textcolor{blue}{b} \\ \textcolor{red}{b-1} & \end{array} 
  \rightarrow
  \begin{array}{c} \textcolor{blue}{b-1} \\ \textcolor{red}{b-1} \end{array}, 
  \qquad \textrm{ or } \qquad 
  \begin{array}{c} \textcolor{blue}{b} \\ \textcolor{red}{b} \end{array}
  \rightarrow
  \begin{array}{cc} & \textcolor{blue}{b} \\ \textcolor{red}{b-1} & \end{array}. 
\end{align}
In the former case, if there was a $\textcolor{blue}{b-1}$ before the move, 
the pair would have been 
$\begin{array}{c} \textcolor{blue}{b-1} \\ \textcolor{red}{b-1} \end{array}$
in the first place.  
Thus, it is not possible to have $\textcolor{blue}{b-1}$ before the move.  
In the latter case, if there is a singleton $(\textcolor{blue}{b-1})$, 
then the move is 
\begin{align}
\nonumber 
  \begin{array}{cc} & \textcolor{blue}{b} \\ 
    (\textcolor{blue}{b-1} ) & \textcolor{red}{b} \end{array}
  \rightarrow
  \begin{array}{cc} (\textcolor{blue}{b-1}) &  \\ 
    (\textcolor{red}{b-1} ) & \textcolor{blue}{b} \end{array}.  
\end{align}
In either case, the weight of the partition decreases by 1, 
and the number of pairs and the number of singletons do not change.  
Backward move on a singleton is $\textcolor{blue}{b}$ 
becoming $(\textcolor{blue}{b-1})$.  
Since the singletons are ahead of the pairs in the base partition, 
we do not need to perform any adjustments.  
A backward move on a singleton decreases the weight of the partition by 1, 
and the number of pairs and the number of singletons do not change.  

For the forward moves on pairs, we have three cases to consider.  
\begin{align}
\nonumber 
  \begin{array}{cc} & \textcolor{blue}{b} \\ \textcolor{red}{b-1} & \end{array} 
  \rightarrow
  \begin{array}{c} \textcolor{blue}{b} \\ \textcolor{red}{b} \end{array}, 
  \qquad \qquad 
  \begin{array}{cc} \textcolor{blue}{b} &  \\ 
    \textcolor{red}{b} & (\textrm{parts } > b+1) \end{array}
  \rightarrow
  \begin{array}{ccc} & (\textcolor{blue}{b+1}) & \\ 
    \textcolor{red}{b} & & (\textrm{parts } > b+1) \end{array}, 
\end{align}
\begin{align}
\nonumber 
  \textrm{ or } \quad 
  \begin{array}{cc} \textcolor{blue}{b} &  \\ 
    \textcolor{red}{b} & (\textcolor{blue}{b+1}) \end{array}
  \rightarrow
  \begin{array}{cc} & (\textcolor{blue}{b+1})  \\ 
    \textcolor{blue}{b} & (\textcolor{red}{b+1}) \end{array}.  
\end{align}
In each case, the weight of the partition increases by 1, 
and the number of pairs and singletons remain constant.  
A forward move on a singleton is pushing $\textcolor{blue}{b}$ 
as $(\textcolor{blue}{b+1})$.  
Since we are pulling singletons last, 
we are pushing singletons first.  
In the base partition, the singletons are already ahead of the pairs, 
so we do not have to perform any adjustments.  
A forward move on a singleton increases the weight of the partition by 1, 
and the number of pairs and the number of singletons do not change.  

A slight variation on the parameters in Theorem \ref{thmYalcinGenFunc} is 
$( B_{11}, B_{22}, B_{12}, D_1, D_2, K_1, K_2, \gamma )$ 
$=( 2, 1, 1, 1, 1, 1, 1, 1 )$.  
In particular, $D_1$ is changed from 2 to 1.  
For the series 
{\allowdisplaybreaks
\begin{align}
\nonumber 
 S_{0,0}(x) & = \sum_{m, n \geq 0} \frac{ q^{ 2 \binom{m+1}{2} + \binom{n+1}{2} + mn }
      x^{m+n} }{ (q; q)_m (q; q)_n } \textrm{ and } \\ 
\nonumber 
 S_{1,0}(x) & = \sum_{m, n \geq 0} \frac{ q^{ 2 \binom{m+1}{2} + \binom{n+1}{2} + mn + m }
      x^{m+n} }{ (q; q)_m (q; q)_n }, 
\end{align}}
the algorithm in Section \ref{secSrchAlgo} gives 
{\allowdisplaybreaks
\begin{align}
\nonumber 
 S_{0,0}(x) & = S_{0,0}(xq) + (xq + xq^2) S_{1,0}(xq) \\ 
\nonumber 
 S_{1,0}(x) & = S_{0,0}(xq) + xq S_{1,0}(xq).  
\end{align}}
Like in Theorem \ref{thmYalcinGenFunc}, 
it is not possible to stick to ordinary partitions, 
so we use colors to interpret the generated partitions.  
The following theorem is proven using structural induction, as well~\cite{Yalcin_thesis}.  

\begin{theorem}
\label{thmYalcinGenFunc2}
 Let $t_1(m, n)$ be the number of partitions of $n$ into $m$ bicolored parts such that 
 \begin{enumerate}[1.]
  \item Both blue and red parts appear at most once.  
  \item $\textcolor{red}{1}$ does not appear as a part.  
  \item For any $\textcolor{blue}{b}$, we cannot have $\textcolor{red}{(b+1)}$ 
    or $\textcolor{red}{(b+2)}$.  
  \item In any sublist of successive parts ending with $\textcolor{blue}{b}$, 
    we cannot have differences $[1,2^{\ast}]$ with successive red parts.  
    Thus, we cannot have $\textcolor{red}{(b-1)} + \textcolor{blue}{b}$, 
    or $\textcolor{red}{(b-3)}$ $+ \textcolor{red}{(b-2)}$ $+ \textcolor{blue}{b}$, 
    or $\textcolor{red}{(b-5)}$ $+ \textcolor{red}{(b-4)}$ 
    $+ \textcolor{red}{(b-2)}$ $+ \textcolor{blue}{b}$.  
 \end{enumerate}
 Let $t_2(m, n)$ be the number of partitions of $n$ into $m$ bicolored parts
 enumerated by $t_1(m, n)$ in which no $\textcolor{red}{2}$ appears.  
\end{theorem}

The four partitions of $n = 4$ satisfying the conditions stipulated by $t_1(m, n)$ are 
\begin{align}
\nonumber 
 \textcolor{blue}{4}, \qquad
 \textcolor{red}{4}, \qquad
 \textcolor{blue}{1} + \textcolor{blue}{3}, \qquad
 \textcolor{red}{2} + \textcolor{blue}{2}. 
\end{align}

Bringing together \eqref{eqYalcin3Regular}, 
Theorem \ref{thmYalcinGenFunc} and \ref{thmYalcinGenFunc2}, 
we prove Theorem \ref{thmYalcinMain}.

\section{Alternating series and inclusion-exclusion}
\label{secAltSeries}

After the discussion in Sections \ref{secSrchAlgo}-\ref{secNarrowingSrchSpc},
allowing the series \eqref{srMainDouble} to be
\begin{align}
\label{srMainDoubleAlt}
 S_{C_1, C_2}(x) = S_{C_1, C_2}(x; q) = \sum_{m, n \geq 0} \epsilon_1^m \epsilon_2^n
  \frac{ q^{ B_{11} \binom{m+1}{2} + B_{22} \binom{n+1}{2} + B_{12}mn + C_1 m + C_2 n }
    x^{ D_1 m + D_2 n } }
   { (q^{K_1}; q^{K_1})_m (q^{K_2}; q^{K_2})_n }
\end{align}
for $\epsilon_1 = \pm 1$, and for $\epsilon_2 = \pm 1$ comes almost for free.
The functional equations are
\begin{align}
\nonumber
 S_{C_1, C_2}(x) - S_{C_1 + K_1, C_2}(x)
 & = \epsilon_1 x^{D_1} q^{ B_{11} + C_1 }
  S_{C_1 + B_{11} - \gamma D_1, C_2 + B_{12} - \gamma D_2 }(x q^{\gamma}) \\
\nonumber
 S_{C_1, C_2}(x) - S_{C_1 , C_2 + K_2}(x)
 & = \epsilon_2 x^{D_2} q^{ B_{22} + C_2 }
  S_{C_1 + B_{12} - \gamma D_1, C_2 + B_{22} - \gamma D_2 }(x q^{\gamma}) \\
\label{fneqPrimaryAlt}
  S_{C_1, C_2}(x) & = S_{C_1 - \gamma D_1 , C_2 - \gamma D_1}(x q^{\gamma}),
\end{align}
accordingly.
The indicated double indices are still the same,
so the auxiliary linear system size estimation detailed in Section \ref{secNarrowingSrchSpc}
stays the same.
Of course, positive series and alternating series have to satisfy different systems of functional equations.
However, we make another deviation from the procedure we followed before.
We will not try and give a partition generating function interpretation
via the system of functional equations we obtain through the algorithm.
This is mainly because systems of functional equations
which evidently imply non-negative coefficients are rarer.
And, it is more difficult to give a signed partition interpretation
via a system of functional equations.
Another reason is that,
once we have an alleged series-product identity,
inclusion-exclusion on a signed partition generating function interpretation
of \eqref{srMainDoubleAlt} is much easier.
It has to be said that many of the identities we came by
were variations on the theme of
\begin{align}
\nonumber
  \sum_{m = 0}^N \begin{bmatrix} N \\ m \end{bmatrix} (-1)^m q^{ \binom{m}{2} } 
  = \begin{cases} 1 & \textrm{ for } N = 0 \\ 0 & \textrm{ for } N > 0 \end{cases}, 
\end{align}
by the $q-$binomial theorem~\cite{A_thebluebook}.  
This can easily be discerned as an inner single sum in the double sum \eqref{srMainDoubleAlt}.
We discarded them.
We prove two of the more interesting examples here, 
namely Theorems \ref{thmAltRRG3} and \ref{thmAltAlladiGordonDilation}.

The search algorithm gave the output in Theorems \ref{thmAltRRG3} and \ref{thmAltAlladiGordonDilation}, 
naturally, as $=$ replaced by $\stackrel{?}{=}$.  

\begin{proof}[proof of Theorem \ref{thmAltRRG3}]
 We prove the third identity.  
 Proofs of others are similar.  
 We rewrite the series as 
 \begin{align}
 \nonumber 
  \sum_{m, n \geq 0} q^{ 2 \binom{m+n+1}{2} } \; 
   \frac{ q^m }{ (q; q)_m } \; \frac{ (-1)^n q^n }{ (q^2; q^2)_n }
 \end{align}
 The only fact we need concerning the first factor $q^{ 2 \binom{m+n+1}{2} }$ 
 in the general term at the moment is that it is a function of $m+n$.  
 The second factor $\frac{ q^m }{ (q; q)_m }$ in the general term 
 generates partitions into $m$ parts, 
 and the third one $\frac{ (-1)^n q^n }{ (q^2; q^2)_n }$ 
 generates partitions into $n$ odd parts, each with weight $(-1)$.  
 Call instances of those partitions $\mu$ and $\nu$, respectively.  
 
 Choose and fix a pair $(\mu, \nu)$ of partitions described above.  
 Also, choose and fix an odd positive integer $(2j-1)$.  
 Let $\mu$ and $\nu$ share a total of $f_{(2j-1)}$ $(2j-1)$'s between them.  
 Then, we can take a $(2j-1)$ from $\mu$, append it to $\nu$ with weight $(-1)$.  
 Or, we can take a $(2j-1)$ from $\nu$, take off its sign, and append it to $\mu$.  
 This creates $f_{(2j-1)}+1$ sister pairs of partitions $(\mu, \nu)$.  
 The sum of their number of parts $m+n$ are the same, 
 and they are counted 
 \begin{align}
 \nonumber 
  \epsilon \left( 1 - 1 + 1 - \cdots + (-1)^{f_{(2j-1)}} \right) 
  = \begin{cases} \epsilon & \textrm{ if } f_{(2j-1)} \textrm{ is even, } \\ 
   0 & \textrm{ if } f_{(2j-1)} \textrm{ is odd } \end{cases}
 \end{align} 
 times.  
 Here, $\epsilon$ is the product of $(-1)$'s 
 coming from parts of $\nu$ other than $(2j-1)$.  
 
 We deduce that the survivors in this inclusion-exclusion 
 are those $(\mu, \nu)$ in which the total count of each odd part is even.  
 Once we agree that $\nu$ can have all of those odd parts, 
 then the survivors are pairs of partitions $(\mu, \nu)$ 
 such that $\mu$ is a partition into $m$ even parts, 
 and $\nu$ is a partition into odd parts appearing evenly~\cite{A_parity}.  
 Each such pair is weighted by $(+1)$ now.  
 We can update the number of parts in $\nu$ as $2n$, and arrive at 
 \begin{align}
 \nonumber 
  \sum_{m, n \geq 0} q^{ 2 \binom{m+n+1}{2} } \; 
   \frac{ q^m }{ (q; q)_m } \; \frac{ (-1)^n q^n }{ (q^2; q^2)_n }
  & = \sum_{m, n \geq 0} q^{ 2 \binom{m+2n+1}{2} } \; 
   \frac{ q^{2m} }{ (q^2; q^2)_m } \; \frac{ q^{2n} }{ (q^4; q^4)_n } \\ 
 \nonumber 
  & = \sum_{m, n \geq 0} 
   \frac{ q^{ 8 \binom{m+1}{2} + 2 \binom{n+1}{2} + 4mn + 2n } }{ (q^4; q^4)_m (q^2; q^2)_n }.  
 \end{align}
 In the last series, we swapped $m \leftrightarrow n$.  
 The last series is the $q \leftarrow q^2$ dilation of 
 the first formula in~\cite[Theorem 1.6]{Yalcin}.  
 This concludes the proof.  
\end{proof}

\begin{proof}[proof of Theorem \ref{thmAltAlladiGordonDilation}]
 The proof is reminiscent of the proof of Theorem \ref{thmAltRRG3}.  
 We rewrite the series in the theorem as 
 \begin{align}
 \nonumber 
  \sum_{m, n \geq 0} q^{ 6 \binom{m+n}{2} } \; 
   \frac{ q^{ 3 \binom{m+1}{2} } }{ (q^3; q^3)_m } \; \frac{ (-1)^n q^n }{ (q^2; q^2)_n }.  
 \end{align}
 The only quality we need of the first factor $q^{ 6 \binom{m+n+1}{2} }$ 
 of the general term is that it is a function of $m+n$.  
 The second factor $\frac{ q^{ 3 \binom{m+1}{2} } }{ (q^3; q^3)_m }$ 
 of the general term generates partitions into $m$ distinct multiples of three, 
 and the third one $\frac{ (-1)^n q^n }{ (q^2; q^2)_n }$ 
 generates partitions into $n$ odd parts, each weighted by $(-1)$.  
 Call those partitions $\mu$ and $\nu$, respectively.  
 
 Choose and fix a pair $(\mu, \nu)$ of partitions as described above.  
 If $\mu$ has an odd multiple of three, say $(6j-3)$, 
 then we can take it, endow it with a $(-1)$, and append it to $\nu$.  
 This gives $(\mu, \nu)$ a sibling pair of partitions 
 with the same sum and total number of parts, generated with with opposite sign.  
 Unlike in the proof of Theorem \ref{thmAltRRG3}, 
 this is doable only once, because $\mu$ has distinct parts.  
 If $\mu$ has no odd parts, but $\nu$ has a part divisible by three, 
 say $(6j-3)$, then we can take it from $\nu$, 
 strip off the weight $(-1)$, and append it to $\mu$.  
 With other things being equal, 
 this causes another cancellation of weights $+1 - 1 = 0$.  
 
 After these cancellations, 
 the survivor pairs are $(\mu, \nu)$ such that 
 $\mu$ is a partition into distinct multiples of six; 
 and $\nu$ is a partition into odd parts that are not multiples of three, 
 each of which is weighted by $(-1)$.  
 In other words, parts of $\nu$ are congruent to $\pm 1 \pmod{6}$.  
 So, it is convenient to separate parts of $\nu$.  
 Let's say it has $n$ parts that are congruent to $1 \pmod{6}$, 
 and $r$ parts that are are congruent to $-1 \pmod{6}$.  
 We now have 
 \begin{align}
 \nonumber 
  \sum_{m, n \geq 0} q^{ 6 \binom{m+n}{2} } \; 
   \frac{ q^{ 3 \binom{m+1}{2} } }{ (q^3; q^3)_m } \; \frac{ (-1)^n q^n }{ (q^2; q^2)_n } 
  & = \sum_{m, n, r \geq 0} q^{ 6 \binom{m+n+r}{2} } \; 
   \frac{ q^{ 6 \binom{m+1}{2} } }{ (q^6; q^6)_m } \; \frac{ (-1)^n q^n }{ (q^6; q^6)_n }
   \; \frac{ (-1)^r q^{5r} }{ (q^6; q^6)_r } \\ 
  \nonumber 
  & = \sum_{m, n, r \geq 0} q^{  } \; 
   \frac{ q^{ 6 \binom{m+n+1}{2} } q^{ 6 \binom{m}{2} } 
      \left( (-q)(-q^5) \right)^m \left( -q \right)^n \left( -q^5 \right)^r }
    { (q^6; q^6)_m (q^6; q^6)_n (q^6; q^6)_r }.  
 \end{align}
 The last series is the $q \leftarrow q^6$ dilation of (2.11) in~\cite{Alladi_Gordon}, 
 with $A = -q$ and $B = -q^5$.  The claim follows.  
\end{proof}

\section{Future work}
\label{secFutureWork}

More extensive search even for double series is needed, 
in which larger bounds for the parameters, 
or larger auxiliary linear systems are taken into account.  
A shortcut for parameters that will likely work is 
trying the Euler Algorithm after substituting $x = q$, $q^2$, $q^3$, \ldots.  
This is twisting the parameters $C_1$ and $C_2$ without running the algorithm all over again.  
An immediate continuation is to rewrite the code to treat $k$-fold sums in general.  
This is ongoing work.  

Another possibility is to use real or even complex characters~\cite{Apostol} instead of 
powers of $(-1)$ like in Section \ref{secAltSeries}.  
Such series previously appeared in~\cite{ABM15, Capp_th, Sills04}.  
This requires adaptation of the primary $q$-contiguous equations 
\eqref{fneqPrimary} or \eqref{fneqPrimaryAlt}.  

Cylindric partition generating functions such as in~\cite{CDU22, CW19} can be considered.  
However; not only the primary $q$-contiguous equations 
\eqref{fneqPrimary} need to be updated, 
but also further optimization is needed.  
This is because the number of series is larger in each such system.  

Checking unique determination of series in functional equations 
by hand takes time.  
The consistency of the systems and the series therein 
with the initial values $S_{\cdot, \cdot}(0) = 1$ 
takes even longer to verify by hand.  
If it is possible to automate these, 
it will save precious time spent by checking redundant systems.  

Decision algorithms such as in~\cite{PWZ96} 
or structure theorems are desired.  
This requires a different viewpoint, 
as it is not possible with the methods used here.  

Last, but not least, a bijective proof of Theorem \ref{thmYalcinMain} 
would be nice.  

\section*{Acknowledgements}

Some parts of this note are from the first author's PhD dissertation.  
We thank Peter Paule, Shashank Kanade and Matthew Russell 
for the useful discussions during the preparation of this manuscript.  
We also thank the anonymous referee 
for carefully reading the manuscript, suggestions for future work, 
and crucial references~\cite{Chern, KR_20}.  


\bibliographystyle{amsplain}

\end{document}